%% file: iop_main.tex
\documentclass{iopjournal}
\usepackage{amsmath}
\usepackage{amsfonts} 
\usepackage{subcaption}  
\usepackage{comment}
\usepackage{tabularx}
\input{macros}
\begin{document}


\title{A Physics-Informed Neural Network with a Modified Lorentzian Activation for Nonlocal Gradient-Flow Equations in Dynamic Density Functional Theory}

\author{Dimitrios Gourzoulidis$^{1,*}$\orcid{0000-0002-3477-3226}, Soumaya Elkantassi$^2$\orcid{0000-0003-2610-3676} and Serafim Kalliadasis$^{1,*}$\orcid{0000-0001-9858-3504}}

\affil{$^1$Department of Chemical Engineering, Imperial College London, United Kingdom}

\affil{$^2$Department of Operations, University of Lausanne,  Switzerland}

\affil{$^*$Author to whom any correspondence should be addressed.}



\begin{abstract}
We develop a physics-informed neural network (PINN) framework for nonlocal partial differential equations arising in dynamic density functional theory (DDFT). Such equations are challenging for standard PINN methods because they involve nonlinearities, nonlocal interaction terms, and an underlying gradient-flow structure, often leading to slow convergence and difficult optimization. We adapt the PINN methodology to DDFT gradient-flow equations and introduce two computational components: a modified Lorentzian activation function that behaves approximately linearly for small inputs and decays toward zero as the input magnitude increases, and a precomputed discrete operator for evaluating the nonlocal convolution term efficiently during training. The method is tested on four examples in one and two space dimensions. In the first example, the exact stationary solution is known, while in the remaining cases the neural-network approximations are validated against reference solutions computed using continuous and discontinuous Galerkin finite element discretizations. Accuracy and physical consistency are assessed through $L^1$, $L^2$, and $L^\infty$ errors, together with mass conservation and free-energy dissipation. The results show that the proposed activation function accelerates convergence relative to the standard \texttt{tanh} function, while the overall framework maintains good agreement with the reference solutions and captures the expected gradient-flow behaviour. These findings demonstrate the potential of the proposed PINN framework for solving nonlocal gradient-flow equations arising in DDFT.
\end{abstract}

\vspace{0.5em}
\noindent
{\large\bfseries Nomenclature}

\vspace{0.5em}

\noindent
\begin{tabularx}{\textwidth}{@{}>{$}l<{$}cX@{}}
\Omega \subset \mathbb{R}^d & = & spatial domain, \\
\partial\Omega & = & boundary of the domain, \\
\mathbf{n} & = & outward unit normal on $\partial\Omega$, \\
\rho(\mathbf{x},t) & = & density, \\
\rho_0(\mathbf{x}) & = & initial density, \\
\mu = \delta\mathcal E/\delta\rho & = & chemical potential, \\
\E[\rho] & = & free-energy functional, \\
\mathcal D(\rho) & = & energy-dissipation functional, \\
H(\rho) & = & local free-energy density, \\
V(\mathbf{x}) & = & prescribed confining potential, \\
W(\mathbf{x}-\mathbf{y}) & = & interaction kernel, \\
W * \rho & = & spatial convolution of $W$ and $\rho$, \\
k_B & = & Boltzmann constant,\\
\theta & = & trainable neural-network parameter vector, \\
\rho_\theta(\mathbf{x},t) & = & neural-network representation of the density, parametrized by the vector $\theta$,\\
(\cdot)^{\T} & = & transpose of a vector or matrix.
\end{tabularx}

\vspace{0.5em}
\noindent
{\large\bfseries Abbreviations}

\vspace{0.5em}

\noindent
\begin{tabularx}{\textwidth}{@{}l c X@{}}
PINN & = & physics-informed neural network, \\
PDE  & = & partial differential equation, \\
DFT  & = & density functional theory, \\
DDFT & = & dynamic density functional theory, \\
FEM    & = & finite element method, \\
CG-FEM & = & continuous Galerkin finite element method, \\
DG-FEM & = & discontinuous Galerkin finite element method.
\end{tabularx}

\newpage
\section{Introduction}

Gradient-flow partial differential equations form an important class of models for dissipative systems, including linear and nonlinear diffusion, aggregation, and nonlocal interaction dynamics \cite{JordanKinderlehrerOtto1998,Villani2003,CarrilloMcCannVillani2003}. Within this framework, dynamic density functional theory (DDFT) provides a continuum description for the evolution of interacting particle densities through a free-energy-driven formulation rooted in classical density functional theory \cite{Evans1979,MarconiTarazona1999,ArcherRauscher2004,GoddardNoldSavvaPavliotisKalliadasis2012,GoddardPavliotisKalliadasis2012}. DDFT and related density-functional approaches have been used to study a broad range of physical phenomena, including colloidal fluids \cite{MarconiTarazona1999,ArcherRauscher2004,GoddardNoldSavvaPavliotisKalliadasis2012}, phase separation \cite{ArcherEvans2004}, wetting and interfacial dynamics \cite{NoldSibleyGoddardKalliadasis2014,YatsyshinSavvaKalliadasis2015}, confined systems \cite{GoddardNoldKalliadasis2016}, and nucleation processes \cite{DuranOlivenciaYatsyshinKalliadasisLutsko2018}.

Despite this appealing variational structure, DDFT-type gradient-flow equations remain numerically challenging. The interplay of nonlinearities, nonlocal interaction operators, and external potentials gives rise to complex multiscale solution behaviour, while the underlying gradient-flow structure imposes fundamental physical constraints, notably mass conservation and free-energy dissipation. Consequently, numerical approximations must achieve not only accuracy but also preserve these essential physical properties. This has led to the development of spectral methods \cite{YatsyshinSavvaKalliadasis2012}, finite-volume methods \cite{CarrilloChertockHuang2015,MendesRussoPerezKalliadasis2021}, continuous Galerkin finite element methods (CG-FEM) \cite{GourzoulidisKalliadasis2026}, and discontinuous Galerkin finite element methods (DG-FEM) \cite{SunCarrilloShu2018} for nonlinear and nonlocal gradient-flow problems, with particular emphasis on preserving the positivity of the numerical density and ensuring discrete energy decay. Although such methods provide reliable numerical approximations, their implementation can become technically demanding in more complex settings, thereby prompting the exploration of alternative computational approaches.

Physics-informed neural networks (PINNs) have emerged as a  used framework for the numerical solution of differential equations by incorporating the governing equations, boundary conditions, and initial data directly into the training objective \cite{LagarisLikasFotiadis1998,SirignanoSpiliopoulos2018,RaissiPerdikarisKarniadakis2019}. Their mesh-free formulation and reliance on automatic differentiation have led to a wide range of extensions, including variational and finite-element-inspired formulations \cite{BerroneCanutoPintore2022}, finite-element-trained neural methods \cite{GrekasMakridakis2025}, time-discrete Runge--Kutta-based approaches \cite{AkrivisMakridakisSmaragdakis2025}, and methods for nonlocal and integro-differential equations \cite{PangDEliaParksKarniadakis2020,YuanNiDengHao2022}. However, to the best of our knowledge, PINNs have not previously been applied to the direct numerical solution of nonlocal DDFT equations. Such problems are particularly challenging because of their nonlinear and nonlocal terms and the need to reproduce the monotonic decay of the free energy.

In this work, we develop a PINN framework for nonlocal PDEs arising in DDFT. In addition to extending PINNs to this class of equations, we introduce two methodological contributions: a modified Lorentzian activation function designed to improve convergence and approximation accuracy, and a precomputed discrete operator for the efficient evaluation of the nonlocal convolution term during training. The proposed method is assessed through four examples in one and two space dimensions: a nonlocal interaction problem with a known stationary solution \cite{CarrilloChertockHuang2015}, a nonlinear diffusion problem in a double-well potential \cite{CarrilloChertockHuang2015}, a nonlinear diffusion problem with nonlocal attractive interactions \cite{MendesRussoPerezKalliadasis2021}, and a two-dimensional droplet coalescence problem \cite{GourzoulidisKalliadasis2026,YatsyshinKalliadasis2021}. For all four examples, the neural-network solutions are compared with reference solutions computed using CG-FEM and DG-FEM discretizations. Accuracy and physical consistency are assessed using $L^1$, $L^2$, and $L^\infty$ errors, together with diagnostic plots of mass conservation and free-energy dissipation. Our numerical results show that the proposed activation function yields faster convergence than the standard \texttt{tanh} activation function while maintaining close agreement with the reference solutions.

The remainder of the paper is organized as follows. Section~\ref{sec:model} introduces the class of nonlocal gradient-flow equations considered in this work and recalls the corresponding free-energy structure. Section~\ref{sec:nn_formulation} presents the PINN formulation and the associated loss functional. Section~\ref{sec:implementation} describes the numerical implementation, including the modified Lorentzian activation function, the discrete residual approximation, and the training procedure. Section~\ref{sec:num} presents the numerical results and compares the proposed approach with CG-FEM and DG-FEM reference solutions. For completeness, Appendix~\ref{AppendixA} describes the construction and use of the precomputed operator for evaluating the nonlocal interaction term, while Appendix~\ref{AppendixB} presents the CG-FEM and DG-FEM discretizations used to compute the reference solutions.

\section{Gradient flow for nonlinear diffusion}\label{sec:model}
Let $\Omega \subset \mathbb{R}^d$ with $d\in\{1,2,3\}$ be a bounded domain with Lipschitz boundary $\partial\Omega$ and outward unit normal $\mathbf{n}$.
We consider a nonlinear diffusion equation in Wasserstein gradient-flow form for a nonnegative density $\rho:\Omega\times[0,T]\to[0,\infty)$.
Given $T>0$, we seek $\rho$ satisfying
\begin{equation}
\partial_t \rho - \nabla \cdot \left( \rho \nabla \mu \right) = 0
\qquad \text{in } \Omega \times (0,T],
\label{eq:main}
\end{equation}
where $\mu$ denotes the variational derivative of the free energy. 
We impose the no-flux boundary condition
\begin{equation}
\mathbf{n}\cdot\left(\rho \nabla \mu\right)=0
\qquad \text{on } \partial\Omega \times (0,T],
\label{eq:noflux}
\end{equation}
which prevents mass transfer across $\partial\Omega$ and ensures the conservation of the total mass $M= \int_\Omega \rho\,d\mathbf{x}$. The initial condition is
\begin{equation}
\rho(\mathbf{x},0)=\rho_0(\mathbf{x})
\qquad \text{in } \Omega,
\label{eq:ic}
\end{equation}
where $\rho_0$ is smooth and strictly positive \cite{Carrillo2019}.

The free energy is
\begin{equation}
\E(\rho)
= \int_{\Omega} H(\rho)\,d\mathbf{x}
+ \int_{\Omega} V(\mathbf{x})\,\rho\,d\mathbf{x}
+ \frac{1}{2}\int_{\Omega} (W*\rho)(\mathbf{x})\,\rho(\mathbf{x})\,d\mathbf{x}.
\label{eq:energy}
\end{equation}
Here $H(\rho)$ is a local entropic term, $V(\mathbf{x})$ is a prescribed confining potential, and $W*\rho$ accounts for nonlocal interactions through the spatial convolution
\begin{equation}
(W*\rho)(\mathbf{x},t)=\int_{\Omega} W(\mathbf{x}-\mathbf{y})\,\rho(\mathbf{y},t)\,d\mathbf{y}.
\label{eq:convolution}
\end{equation}
Assuming sufficient regularity, the variational derivative can be written explicitly as
\begin{equation}
\mu := \dfrac{\delta \E}{\delta \rho}
= H'(\rho) + V + (W*\rho).
\label{eq:variation}
\end{equation}
Substituting \eqref{eq:variation} into \eqref{eq:main} yields the equivalent form
\begin{equation}
\partial_t \rho
= \nabla\cdot\Big(\rho \nabla\big[H'(\rho)+V+W*\rho\big]\Big).
\label{eq:equivalent}
\end{equation}
We assume $H \in C^{1}((0,\infty)) \cap C([0,\infty))$ is convex, $V \in C^{1}(\overline{\Omega})$, and $W \in C^{1}(\mathbb{R}^d)$ is symmetric in the sense that $W(-\mathbf{z})=W(\mathbf{z})$ for all $\mathbf{z}\in\mathbb{R}^d$.
In the flux appearing in \eqref{eq:main}, the factor $\rho$ acts as the mobility in the Wasserstein gradient-flow formulation.
Since this mobility vanishes when $\rho=0$, the equation may be degenerate in low-density regions, which can affect the regularity theory and well-posedness properties of the problem.

For smooth solutions of \eqref{eq:equivalent} satisfying \eqref{eq:noflux} and \eqref{eq:ic}, the free energy dissipates according to
\begin{equation}
\dfrac{d}{dt}\E[\rho(t)] = -\mathcal{D}(\rho(t)),
\qquad
\mathcal{D}(\rho) := \int_\Omega \rho \left| \nabla \dfrac{\delta{\E}}{\delta\rho} \right|^2 \, d\mathbf{x} \geq 0,
\label{eq:energy_dissipation}
\end{equation}
so that $\E(\rho(t))$ is non-increasing in time, while \eqref{eq:noflux} ensures mass conservation.

The nonlinear and nonlocal gradient-flow equations introduced above include models arising in dynamic density functional theory. For suitable choices of the local free-energy density $H$, the external potential $V$, and the interaction kernel $W$, equation~\eqref{eq:equivalent} recovers the standard overdamped DDFT equation for the density of interacting Brownian particles. Such formulations are well established in the DDFT literature; see, for example, \cite{ArcherEvans2004,Evans2016,GoddardNoldSavvaPavliotisKalliadasis2012,GoddardNoldKalliadasis2016}.  In the numerical examples, we consider nonlocal interaction models, nonlinear diffusion in external potentials, aggregation--diffusion equations with attractive interactions, and a DDFT model of droplet coalescence.


\section{Neural network formulation} \label{sec:nn_formulation}
In the PINN framework, we assume that the solution $\rho$ belongs to the parametric family $\mathcal{N}=\{\rho_\theta:\theta\in\Theta\}$ and estimate the parameter vector $\theta$ by minimizing an $L^2$ empirical objective built from the PDE residual together with the initial and boundary conditions \cite{RaissiPerdikarisKarniadakis2019,KarniadakisKevrekidisLuPerdikarisWangYang2021}. From a statistical perspective, this corresponds to empirical risk minimization over the finite-dimensional parameter space $\Theta$.

We take $\rho_\theta$ to be a fully connected feedforward neural network with input dimension $d_1$, hidden widths $d_2,\dots,d_L$, output dimension $d_{L+1}$, and $L$ affine transformations, defined by the composition
\begin{equation}\label{eq:rho-NN-composition}
\rho_\theta(\mathbf{x},t)
= \big(C_L \circ \sigma \circ C_{L-1} \circ \cdots \circ \sigma \circ C_1\big)(\mathbf{x},t),
\end{equation}
where $\sigma$ is an activation function applied componentwise. The maps $C_k:\mathbb{R}^{d_k}\to\mathbb{R}^{d_{k+1}}$ are affine and given by
\begin{equation}\label{eq:rho-affine-layers}
C_k(\mathbf{z}) = W_k \mathbf{z} + \mathbf{b}_k,
\qquad
W_k \in \mathbb{R}^{d_{k+1}\times d_k},\;\; \mathbf{b}_k \in \mathbb{R}^{d_{k+1}},
\end{equation}
so that the parameter vector $\theta=\{W_k,\mathbf{b}_k\}_{k=1}^{L}$ consists of all trainable weight matrices and bias vectors of the network and has dimension
\begin{equation}\label{eq:theta-dim}
\dim(\theta)=\sum_{k=1}^{L} d_{k+1}(d_k+1).
\end{equation}

If the activation function $\sigma$ is smooth, then the neural-network representation inherits the corresponding regularity \cite{AkrivisMakridakisSmaragdakis2025}. In practice, we require $\rho_\theta$ to be differentiable with respect to $\mathbf{x}$ and $t$ so that the PDE residual can be evaluated by automatic differentiation \cite{Baydin2018}.

As in standard PINN formulations, we train $\rho_\theta$ by minimizing a composite loss built from terms associated with the governing equation, the initial data, and the boundary condition. The precise form of the residual and the relative weighting of the loss terms may vary across applications and problem classes \cite{RaissiPerdikarisKarniadakis2019,KarniadakisKevrekidisLuPerdikarisWangYang2021}. For the present model, we adopt the strong-form residual
\begin{equation}
r(\rho)(\mathbf{x},t)
:= \partial_t \rho(\mathbf{x},t)
   - \nabla\!\cdot\!\left(\rho(\mathbf{x},t)\,\nabla \mu(\mathbf{x},t)\right),
\label{eq:strong-residual}
\end{equation}
and define the loss functional by
\begin{equation}
\mathcal{L}(\rho)
:= \lambda_1 \!\int_{0}^{T}\!\!\int_{\Omega} \!\left|r(\rho)(\mathbf{x},t)\right|^{2}\,\mathrm{d}\mathbf{x}\,\mathrm{d}t
 + \lambda_2 \!\int_{\Omega}\! \left|\rho(\mathbf{x},0)-\rho_{0}(\mathbf{x})\right|^{2}\,\mathrm{d}\mathbf{x}
 + \lambda_3 \!\int_{0}^{T}\!\!\int_{\partial\Omega}
 \left|\mathbf n\!\cdot\!\left(\rho\,\nabla \tfrac{\delta\mathcal{E}[\rho]}{\delta \rho}\right)\right|^{2}\,\mathrm{d}S\,\mathrm{d}t .
\label{eq:loss-functional}
\end{equation}
Here, the loss functional $\mathcal{L}$ consists of three weighted contributions. The first term penalizes the PDE residual in the interior of the space--time domain, the second penalizes deviations from the initial condition, and the third penalizes violations of the boundary condition through the normal flux on $\partial\Omega$. The positive coefficients $\lambda_1$, $\lambda_2$, and $\lambda_3$ account for differences in scale among these terms and influence both the convergence of the optimization and the quality of the final approximation.

Having defined the loss functional $\mathcal{L}$, we now consider its minimization over the trial set of neural-network approximations
$$
V_N := \left\{ \rho_\theta : \Omega\times[0,T]\to\mathbb{R} \;\big|\; \theta\in\Theta \right\},
$$
that is,
$$
\min_{\rho \in V_N} \; \mathcal{L}(\rho).
$$
Since the trial set $V_N$ is parameterized by the neural-network parameters, this problem is written as the finite-dimensional optimization problem
$$\min_{\theta\in\Theta}\mathcal{L}(\rho_\theta).$$
This is generally a nonconvex optimization problem in $\theta$. In practice, the problem is solved approximately by first-order stochastic optimization methods, with the aim of identifying a parameter vector $\theta$ for which the PDE residual, the initial-condition residual, and the boundary-condition residual are all simultaneously minimized. Moreover, because the integrals in $\mathcal{L}$ are generally not available in closed form, they are approximated numerically, which yields a computable discrete loss. The resulting discrete formulation is described in Section \ref{sec:implementation}.

\section{Implementation} \label{sec:implementation}

This section presents the architecture and numerical realization of the neural-network formulation introduced in Section~\ref{sec:nn_formulation}. We describe the parameterization, the discrete approximations of the chemical potential and residuals, the resulting discrete loss functional, and the optimization procedure used in the computations.

We use the feedforward architecture introduced in Section~\ref{sec:nn_formulation} to define the raw network output $g_\theta$. The density approximation is then obtained by applying a final output transformation to $g_\theta$, as described below. The hidden layers use the modified Lorentzian activation function
$$
\sigma_\varepsilon(z)
=
z\,\dfrac{\varepsilon^2}{z^2+\varepsilon^2}.
$$
Using the affine maps $C_k$ and the parameter vector $\theta$ defined in Section~\ref{sec:nn_formulation}, the raw network output is
$$
g_\theta(\mathbf{x},t)
=
\left(
C_L \circ \sigma_\varepsilon \circ C_{L-1}
\circ \cdots \circ
\sigma_\varepsilon \circ C_1
\right)(\mathbf{x},t).
$$

\begin{figure}[h]
    \centering
    \includegraphics[width=0.49\textwidth]{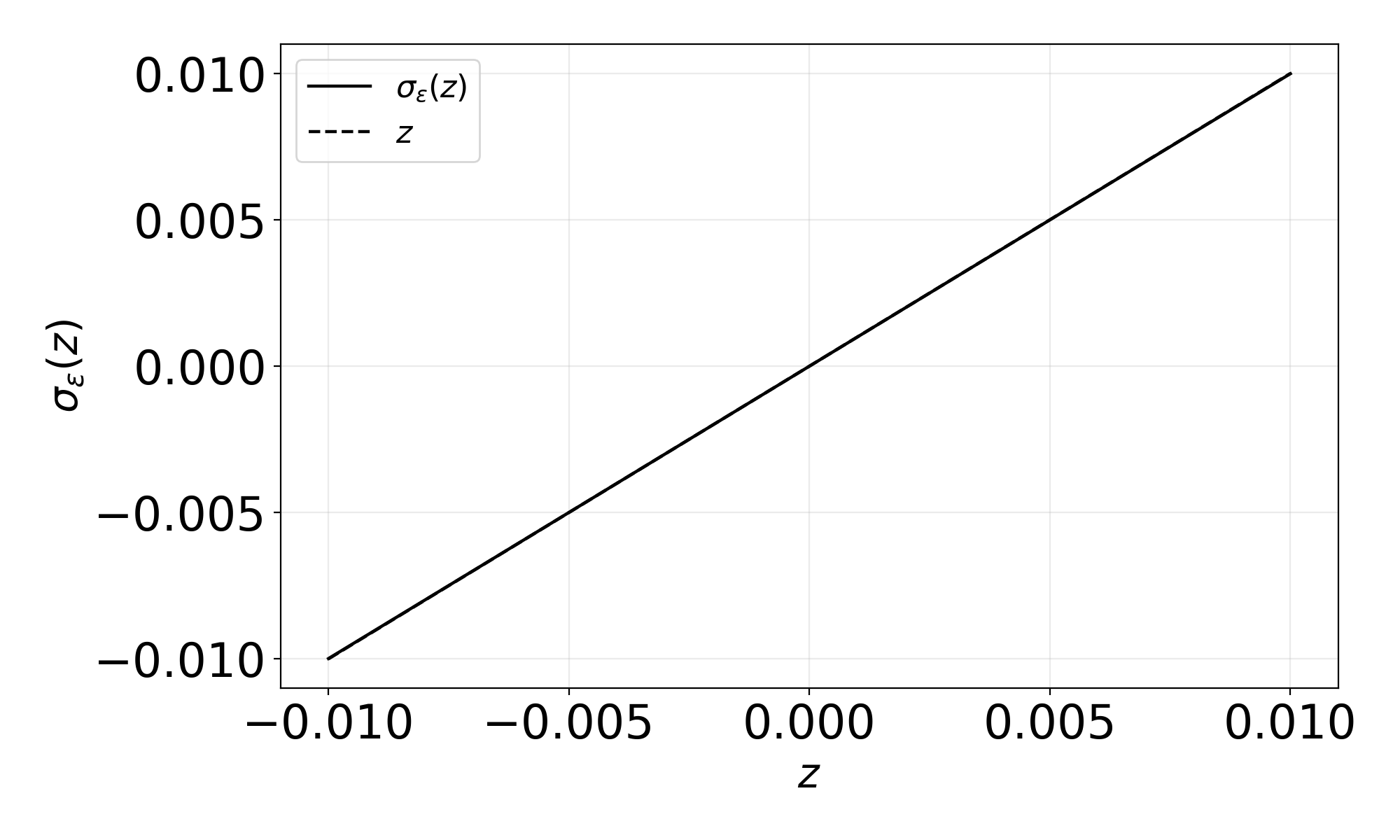}
    \includegraphics[width=0.49\textwidth]{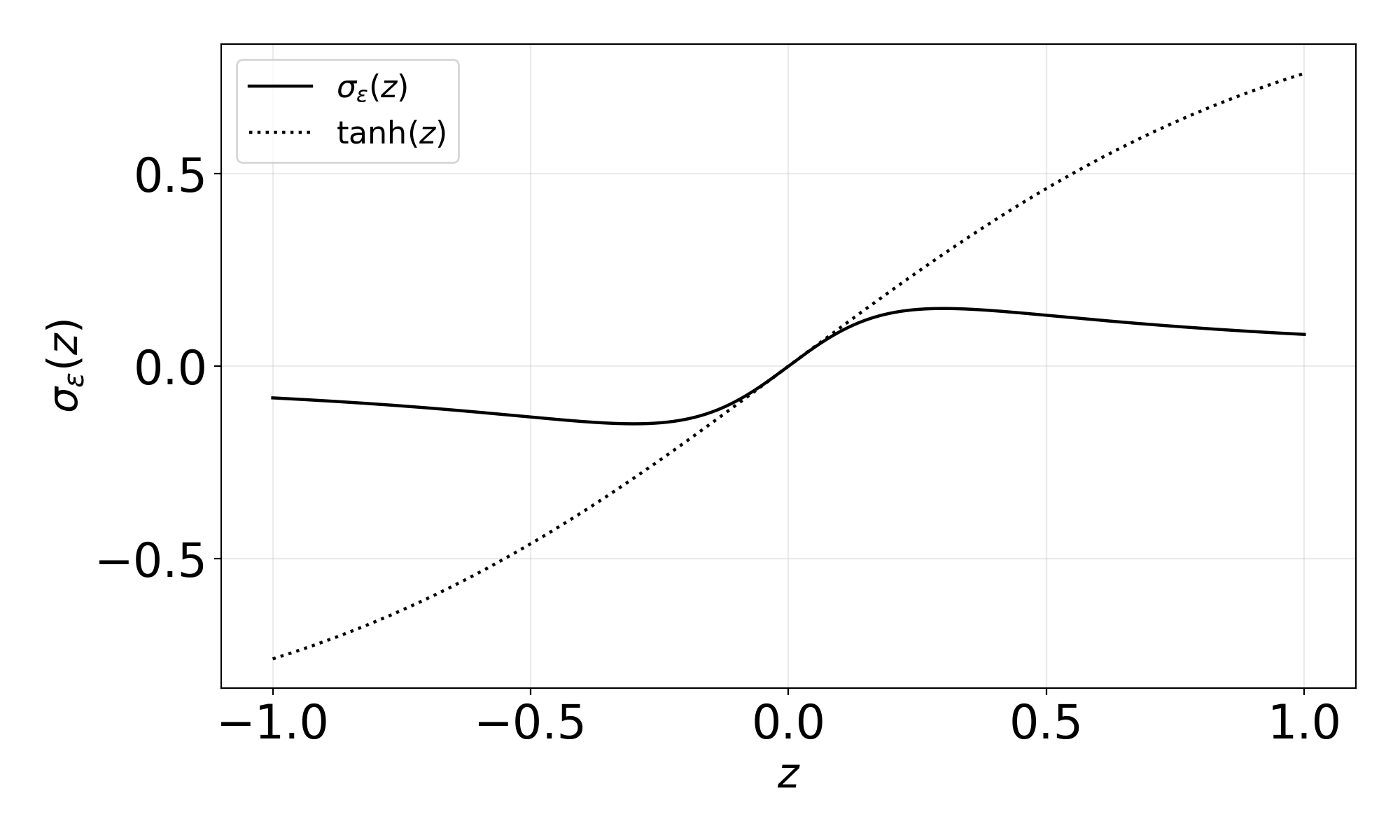}
    \caption{Activation profile of $\sigma_{\varepsilon}(z)$ for $\varepsilon=0.3$. Left: zoom near the origin. Right: comparison with $\tanh(z)$ on a wider interval.}
    \label{fig:activation_sigma}
\end{figure}

Figure~\ref{fig:activation_sigma} illustrates the behaviour of the modified Lorentzian activation near the origin and compares its profile with that of $\tanh(z)$. The activation has two characteristic regimes. For $|z|\ll\varepsilon$, one has $\sigma_\varepsilon(z)\approx z$, so it behaves approximately linearly near the origin. For $|z|\gg\varepsilon$, one has $\sigma_\varepsilon(z)\approx \varepsilon^2/z$, and therefore decays toward zero as the input magnitude increases. Thus, $\sigma_\varepsilon$ combines near-origin linearity with attenuation of large input values. In the numerical experiments, replacing $\tanh$ with $\sigma_\varepsilon$ leads to faster convergence and requires fewer optimization steps to reach a given accuracy, although each iteration has a slightly higher computational cost.

The final density approximation is obtained by applying an output map $\mathcal T$ to the raw output $g_\theta$. In other words,
$$
\rho_\theta(\mathbf{x},t)=\mathcal T\big(g_\theta(\mathbf{x},t)\big),
$$
In practice, the identity map is sufficient in regimes where training remains stable and the learned solution remains nonnegative over the domain. However, during the early stages of training the raw output $g_\theta(\mathbf{x},t)$ may transiently become negative before the parameters approach a physically admissible solution. This may violate the constraint $\rho > 0$ and cause numerical instabilities in the evaluation of the residual. In such cases, we instead use a nonnegative output map, namely the softplus or quadratic transformation, to improve robustness. The output maps considered are
$$
\mathcal T(z)=
\begin{cases}
z, & \text{identity},\\[2pt]
z^2, & \text{quadratic},\\[2pt]
\log(1+e^z), & \text{softplus}.
\end{cases}
$$

To construct discrete approximations of the chemical potential~\eqref{eq:variation}
and the residual~\eqref{eq:strong-residual}, we discretize the spatial domain
$\Omega\subset\mathbb{R}^d$ using a structured Cartesian grid with nodes
$\{\mathbf{x}_i\}_{i=1}^{N}\subset\bar{\Omega}$, where
$\mathbf{x}_i\in\mathbb{R}^d$. We also sample $K$ time collocation points
$\{t_k\}_{k=1}^{K}\subset[0,T]$ independently from the uniform distribution
$\mathcal{U}(0,T)$. For each sampled time $t_k$, the nodal values of the density approximation are
collected in the vector
$$
\boldsymbol{\rho}_\theta(t_k)
=
\left(
\rho_\theta(\mathbf{x}_1,t_k),
\dots,
\rho_\theta(\mathbf{x}_N,t_k)
\right)^{\T}
\in\mathbb{R}^{N}.
$$

The nonlocal interaction term in~\eqref{eq:variation} is evaluated using a
precomputed discrete operator $\Phi\in\mathbb{R}^{N\times N}$, which acts on
$\boldsymbol{\rho}_\theta(t_k)\in\mathbb{R}^{N}$ according to
$$
(W * \rho_\theta)(\mathbf{x}_i,t_k)
\approx
\left(\Phi\,\boldsymbol{\rho}_\theta(t_k)\right)_i,
\qquad i=1,\dots,N.
$$
The construction of $\Phi$ and its relation to the interaction kernel $W$ are described in Appendix~\ref{AppendixA}.

Accordingly, the discrete chemical potential at the grid nodes
$\mathbf{x}_i$ and sampled times $t_k$ is defined by
$$
\mu_\theta(\mathbf{x}_i,t_k)
=
H'\!\left(\rho_\theta(\mathbf{x}_i,t_k)\right)
+
V(\mathbf{x}_i)
+
\left(\Phi\,\boldsymbol{\rho}_\theta(t_k)\right)_i.
$$
Using this discrete chemical potential, the discrete PDE residual at
$\mathbf{x}_i$ and $t_k$ is defined by
$$
r_h(\rho_\theta)(\mathbf{x}_i,t_k)
:=
\partial_t \rho_\theta(\mathbf{x}_i,t_k)
-
\nabla_h\!\cdot\!\left(
\rho_\theta(\mathbf{x}_i,t_k)\,
\nabla_h \mu_\theta(\mathbf{x}_i,t_k)
\right).
$$

The no-flux boundary residual is defined on
$\partial\Omega_h\times\{t_k\}$ by
$$
b_h(\rho_\theta)(\mathbf{x}_i,t_k)
:=
\mathbf{n}(\mathbf{x}_i)\cdot
\mathbf{J}_h(\rho_\theta)(\mathbf{x}_i,t_k),
\qquad
\mathbf{x}_i\in\partial\Omega_h,
$$
where $\mathbf{J}_h(\rho_\theta)$ denotes the discrete approximation of
$\rho_\theta\nabla\mu_\theta$ used in the evaluation of the PDE residual.
Finally, the initial condition is enforced at the spatial grid nodes by penalizing the residual
$$
c_h(\rho_\theta)(\mathbf{x}_i)
:=
\rho_\theta(\mathbf{x}_i,0)-\rho_0(\mathbf{x}_i).
$$

In the discrete residuals above, the temporal derivative
$\partial_t\rho_\theta$ appearing in $r_h$ is computed by automatic
differentiation with respect to $t$, while the spatial derivatives entering
both $r_h$ and the boundary flux $\mathbf{J}_h$ are approximated by finite
differences. More precisely, second--order centred finite differences are used
at interior grid points, while second--order one-sided finite differences are
used at boundary points. The discrete divergence $\nabla_h\!\cdot$ is evaluated
in flux form using the corresponding finite-difference approximations.

The three terms of the loss functional \eqref{eq:loss-functional} are approximated by empirical averages over the spatial grid and sampled times. In practice, we further augment the resulting discrete objective by an additional mass-penalty term in order to improve conservation of the total mass during training. The discrete training objective therefore takes the form
\begin{align*}
\mathcal{L}_h(\theta)
&=
\lambda_1\,
\dfrac{1}{KN}\sum_{k=1}^{K}\sum_{i=1}^{N}
\left|r_h(\rho_\theta)(\mathbf{x}_i,t_k)\right|^2
+
\lambda_2\,
\dfrac{1}{N}\sum_{i=1}^{N}
\left|c_h(\rho_{\theta})(\mathbf{x}_i)\right|^2
\notag\\
&\quad+
\lambda_3\,
\dfrac{1}{K|\partial\Omega_h|}\sum_{k=1}^{K}\sum_{\mathbf{x}_i\in\partial\Omega_h}
\left|b_h(\rho_\theta)(\mathbf{x}_i,t_k)\right|^2
+
\lambda_{4} \dfrac{1}{K} \sum_{k=1}^K
\left| \sum_{i=1}^N \rho_\theta(\mathbf{x}_i,t_k)\,\Delta V - M_0 \right|^2 .
\end{align*}
Here, $|\partial\Omega_h|$ denotes the number of boundary nodes, $\Delta V$ is the grid-cell volume element, and
$ M_0 = \int_\Omega \rho_0\, d\mathbf{x}$
is the prescribed total initial mass.

We minimize $\mathcal{L}_h(\theta)$ using the Adam optimizer. At each iteration, we retain the full Cartesian spatial grid and sample only in time by drawing
$$
\mathbf{t}
=
(t_1,\dots,t_K)^{\T}
\in\mathbb{R}^{K},
\qquad
t_k\overset{\mathrm{i.i.d.}}{\sim}\mathcal{U}(0,T).
$$
The network is then evaluated at all $N$ spatial grid nodes for each sampled time, and the resulting values are arranged in the tensor
$$
R
=
\left[
\rho_\theta(\mathbf{x}_i,t_k)
\right]
\in
\mathbb{R}^{K\times n_1\times\cdots\times n_d},
$$
where $n_j$ denotes the number of grid points in the $j$th spatial direction and
$$
N=\prod_{j=1}^{d}n_j.
$$
This tensor is used to evaluate the PDE, boundary, and mass terms in
$\mathcal{L}_h(\theta)$.

\subsection{Reference finite element solvers}\label{sec:refs}
To assess the neural-network solutions, we compute reference solutions using finite element discretizations implemented in FEniCS~\cite{fenics2015,fenicsbook}. For one-dimensional domains, we use a discontinuous Galerkin discretization with piecewise linear elements on structured meshes~\cite{riviere2008dg,hesthavenwarburton2007dg}. Local conservation is enforced through numerical fluxes at element interfaces, and a minmod slope limiter~\cite{leveque2002finite} is applied after each solve to control spurious oscillations. The nonlocal interaction term is evaluated using a precomputed discrete operator constructed following the same procedure described in Appendix~\ref{AppendixA}, but assembled on the finer mesh used for the finite element computation. This ensures that the PINN and finite element methods use the same discretization principle for the nonlocal term, so that the comparison primarily reflects differences between the solution strategies rather than differences in the approximation of the convolution. Time integration is performed with a semi-implicit Crank--Nicolson scheme, and the nonlinear systems arising at each time step are solved by Newton's method using automatic differentiation within FEniCS/UFL.

For two-dimensional problems, we use a CG-FEM discretization with continuous Lagrange elements on triangular meshes. The nonlocal interaction term, temporal discretization, and nonlinear solution procedure follow the approach described above for the one-dimensional case. The neural-network solutions are assessed against these reference solutions using the error norms reported in Section~\ref{sec:num}. Appendix~\ref{AppendixB} presents the corresponding weak formulations and reference discretizations.

\section{Numerical results} \label{sec:num}

In this section, we present a series of numerical experiments designed to assess the qualitative performance of the proposed method across a range of gradient-flow and DDFT-type problems. The tests include cases with known equilibrium states, which allow direct validation against exact solutions, as well as more challenging interacting-particle configurations. Our goal is to examine whether the method correctly reproduces the expected transient dynamics and preserves the underlying free-energy dissipation mechanism.
The numerical experiments were carried out on a laptop equipped with an NVIDIA RTX 3500 Ada Generation GPU.

\subsection{Example 1: nonlocal interaction test}
We consider a one-dimensional model in which the free energy contains only the nonlocal interaction contribution, with interaction kernel
$$W(x)=\dfrac12 |x|^2-\log |x|,$$
while the internal energy contribution and the external potential are both zero, that is, $$H(\rho)=0, \qquad V(x)=0.$$
The governing equation is therefore
\begin{equation}
\partial_t \rho
=
\partial_x \!\left(
\rho\,\partial_x
\left[
\int_{\Omega}
\left(
\frac{1}{2}|x-y|^2-\log|x-y|
\right)
\rho(y,t)\,dy
\right]
\right),
\qquad
(x,t)\in\Omega\times(0,T].
\label{eq:nonlocal-interaction-test}
\end{equation}

For this choice of interaction kernel, the unique stationary solution of unit mass is known explicitly; see, for instance, \cite{CarrilloChertockHuang2015}, and is given by
$$
\rho_\infty(x)=
\begin{cases}
\displaystyle \dfrac{1}{\pi}\sqrt{2-x^2}, & |x|\le \sqrt{2},\\[4pt]
0, & \text{otherwise}.
\end{cases}
$$
The equilibrium profile \(\rho_\infty\) is \(1/2\)-H\"older continuous. It is the global minimizer of the free-energy functional \eqref{eq:energy}, and solutions of \eqref{eq:main} converge exponentially fast towards this equilibrium; see \cite{CarrilloChertockHuang2015} . In the numerical experiments, we approximate $\rho_\infty$ by solving \eqref{eq:main} up to large times from the Gaussian initial condition
$$
\rho(x,0)=\dfrac{1}{\sqrt{2\pi}}e^{-\frac{x^2}{2}}.
$$

For this experiment, the computations are carried out on the domain $\Omega=[-3,3]$ using a fully connected feedforward network with $4$ hidden layers of width $128$. The loss weights are set to $\lambda_1=1$, $\lambda_2=5$, $\lambda_3=1$, and $\lambda_{4}=0.01$, and the batch size is $100$. The final density is obtained from the raw network output $g_\theta(x,t)$ via the quadratic map. For the modified Lorentzian activation, the parameter $\varepsilon$ is set to $0.05$ The modified Lorentzian model is trained for $10{,}000$ iterations, while the corresponding \texttt{tanh} model is trained for $30{,}000$ iterations.

\begin{figure}[htbp]
  \centering

  \begin{subfigure}{0.48\textwidth}
    \centering
    \includegraphics[width=\textwidth]{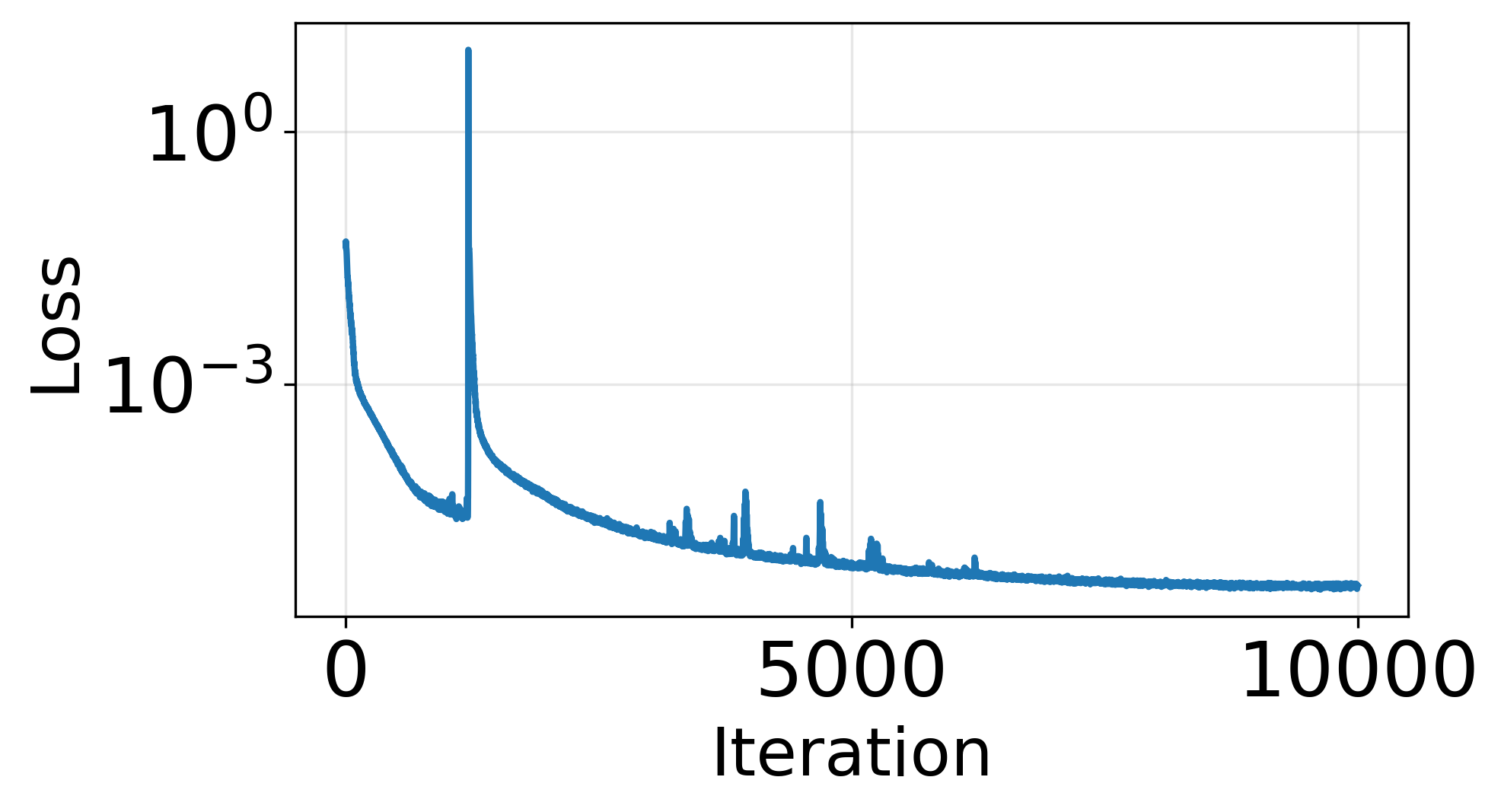}
    \caption{Training loss for the modified Lorentzian activation.}
    \label{fig:test1-loss-lorentz}
  \end{subfigure}
  \hfill
  \begin{subfigure}{0.48\textwidth}
    \centering
    \includegraphics[width=\textwidth]{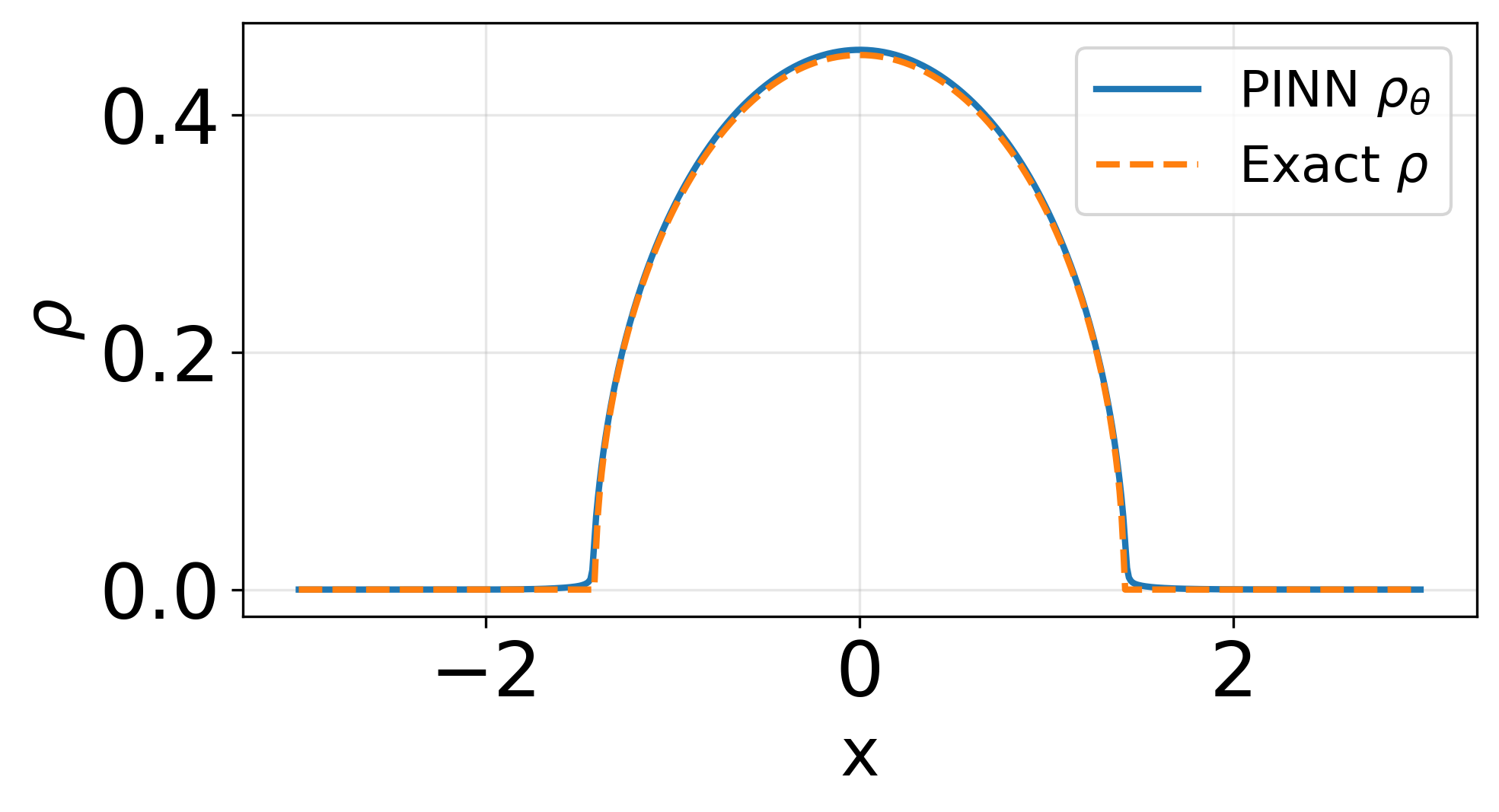}
    \caption{Modified Lorentzian PINN vs exact solution.}
    \label{fig:test1-solution-lorentz}
  \end{subfigure}

  \vspace{0.4cm}

  \begin{subfigure}{0.48\textwidth}
    \centering
    \includegraphics[width=\textwidth]{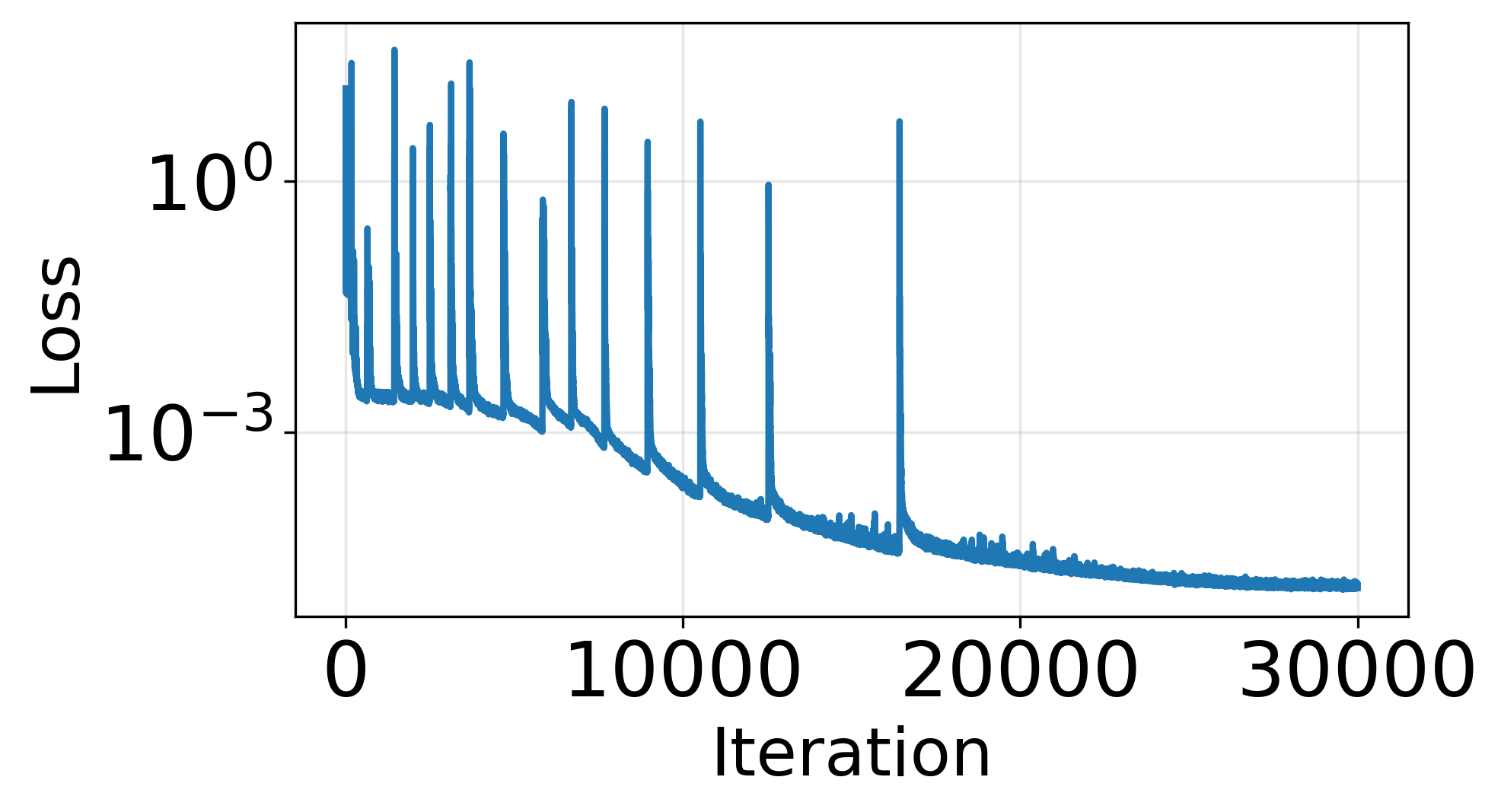}
    \caption{Training loss for the standard \texttt{tanh} activation.}
    \label{fig:test1-loss-tanh}
  \end{subfigure}
  \hfill
  \begin{subfigure}{0.48\textwidth}
    \centering
    \includegraphics[width=\textwidth]{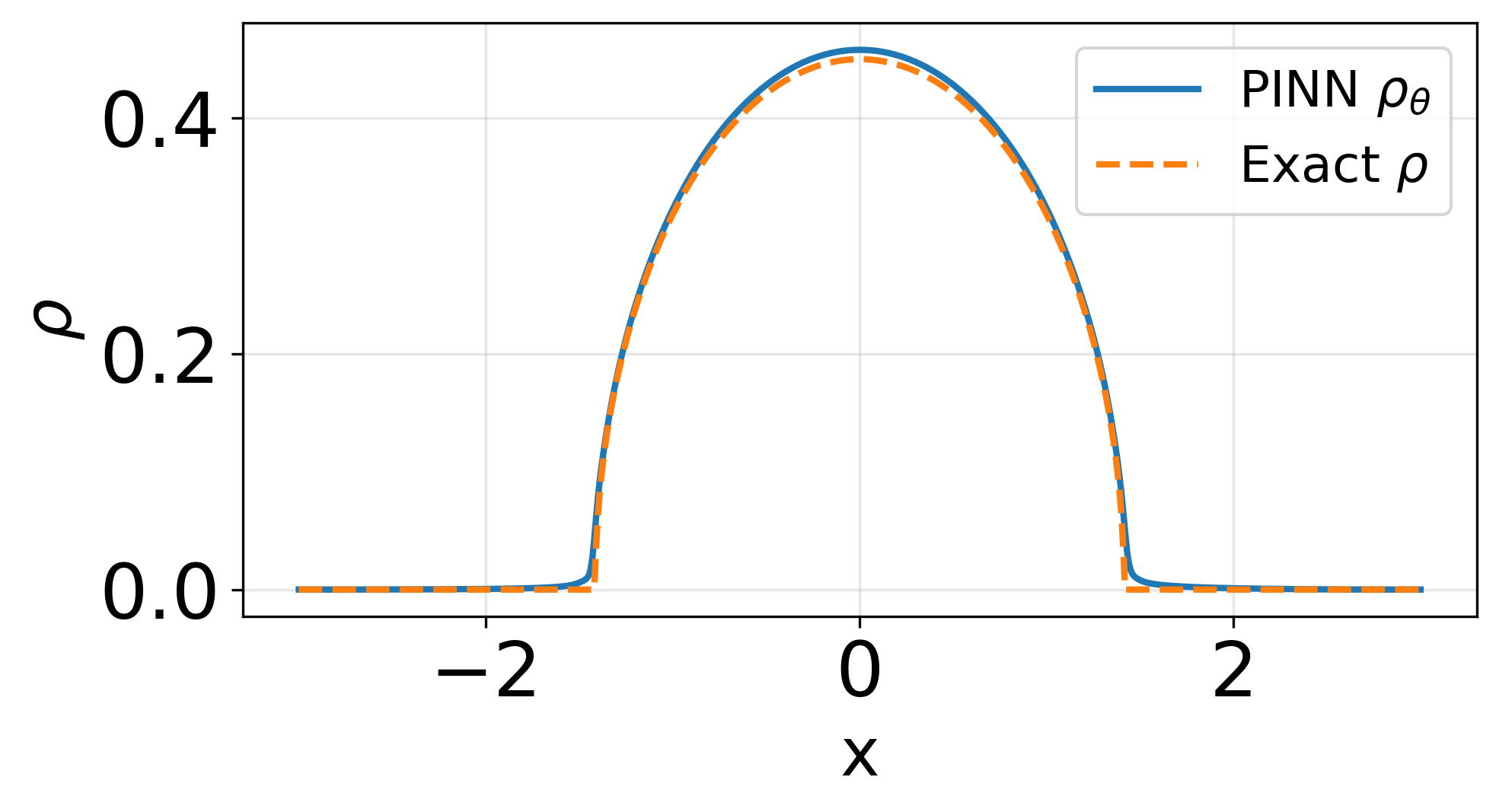}
    \caption{\texttt{tanh} PINN vs exact solution.}
    \label{fig:test1-solution-tanh}
  \end{subfigure}
\caption{Comparison of the modified Lorentzian and standard \texttt{tanh} activations. Panels (a) and (c) show the loss curves during training, while panels (b) and (d) compare the final steady-state approximations with the exact stationary solution.}
  \label{fig:test1-pinn-activation-comparison}
\end{figure}

The left panels of Figure~\ref{fig:test1-pinn-activation-comparison} show the training-loss curves, while the right panels compare the final steady-state approximations obtained using the modified Lorentzian and standard \texttt{tanh} activation functions with the exact stationary solution. All computations are performed on a uniform grid with spacing $\Delta x=0.01$.  Both activation functions produce numerical approximations in good agreement with the exact stationary solution. However, the modified Lorentzian activation yields a more accurate final approximation, with the improvement being particularly visible near $x=0$. Its $L^1$ error is of order $10^{-2}$ and is approximately two times smaller than the $L^1$ error obtained with the \texttt{tanh} activation, while its $L^\infty$ error is also lower.

The main difference, however, lies in the computational cost. The \texttt{tanh} model reaches a loss of order $10^{-5}$ only after $30{,}000$ training iterations, whereas the modified Lorentzian activation reaches the same level after about $3{,}000$ iterations. In wall-clock time, \texttt{tanh} requires about $532\,\mathrm{s}$ to reach this tolerance, whereas the modified Lorentzian activation reaches it in about $236\,\mathrm{s}$. This reduces the training time by more than half. Moreover, after $10{,}000$ iterations, the modified Lorentzian activation reaches a lower loss, of order $10^{-6}$.

Although each iteration with the modified Lorentzian activation is more expensive, about $3.5$ times the cost of one \texttt{tanh} iteration on our hardware, its much faster convergence more than compensates for the larger per-iteration cost. Overall, for this test problem, the modified Lorentzian activation provides both improved optimization efficiency and improved final approximation accuracy relative to the standard \texttt{tanh} baseline.

\subsection{Example 2: nonlinear diffusion with a double-well potential}

We next consider a one-dimensional nonlinear diffusion problem with a double-well external potential, following the example introduced in~\cite{CarrilloChertockHuang2015}. The governing equation is
\begin{equation*}
\partial_t \rho
=
\partial_x\!\left(
\rho\,\partial_x\left(\nu \rho^{m-1} + V(x)\right)
\right),
\qquad
V(x)=\frac{x^4}{4}-\frac{x^2}{2},
\end{equation*}
defined on the domain $\Omega=[-3,3]$ up to the final time $T=4$. 
We set $\nu=1$ and $m=2$, and take a Gaussian density as the initial condition 
\begin{equation*}
\rho_0(x)
=
\dfrac{M}{\sqrt{2\pi\sigma^2}}
\exp\!\left(
-\dfrac{x^2}{2\sigma^2}
\right),
\end{equation*}
with $M=0.1$ and variance $\sigma^2=0.2$. The corresponding free energy is
$$
\E[\rho]=\int_{\Omega}\left(\frac{1}{2}\rho^2+V(x)\rho\right)\,dx,
$$
and $\rho_\infty$ denotes the associated equilibrium density.

This test is of interest because the double-well potential drives an initially centered density towards a symmetric bimodal equilibrium localized near the minima of the external potential.

For this experiment, the loss weights are set to $\lambda_1=1$, $\lambda_2=1$, $\lambda_3=1$, and $\lambda_{4}=0.01$. The neural-network approximation uses a fully connected feedforward architecture with four hidden layers of width $128$, as described in Section~\ref{sec:implementation}. The modified Lorentzian activation is used with parameter $\varepsilon=0.05$. The final density is obtained from the raw network output via the identity map. Training is performed for $50{,}000$ iterations. Here, $\rho_{\mathrm{ref}}:=\rho_{\mathrm{DG}}$ denotes the reference solution.

The neural-network approximation is compared with the reference solution over the full time horizon. The resulting space--time errors are small, indicating good agreement throughout the simulated interval. In particular, the $L^2(0,T;L^2(\Omega))$ error, defined by integration in both space and time, and the corresponding $L^1(0,T;L^1(\Omega))$ error are given by
\begin{align*}
||\rho_{\mathrm{ref}}(x,t)-\rho_{\theta}(x,t)||_{L^2(0,T;L^2(\Omega))}&= 4.42\times 10^{-3},\\
||\rho_{\mathrm{ref}}(x,t)-\rho_{\theta}(x,t)||_{L^1(0,T;L^1(\Omega))}&= 1.13\times 10^{-2}
\end{align*}

\begin{figure}[h]
  \centering

  \begin{subfigure}{0.48\textwidth}
    \centering
    \includegraphics[width=\textwidth]{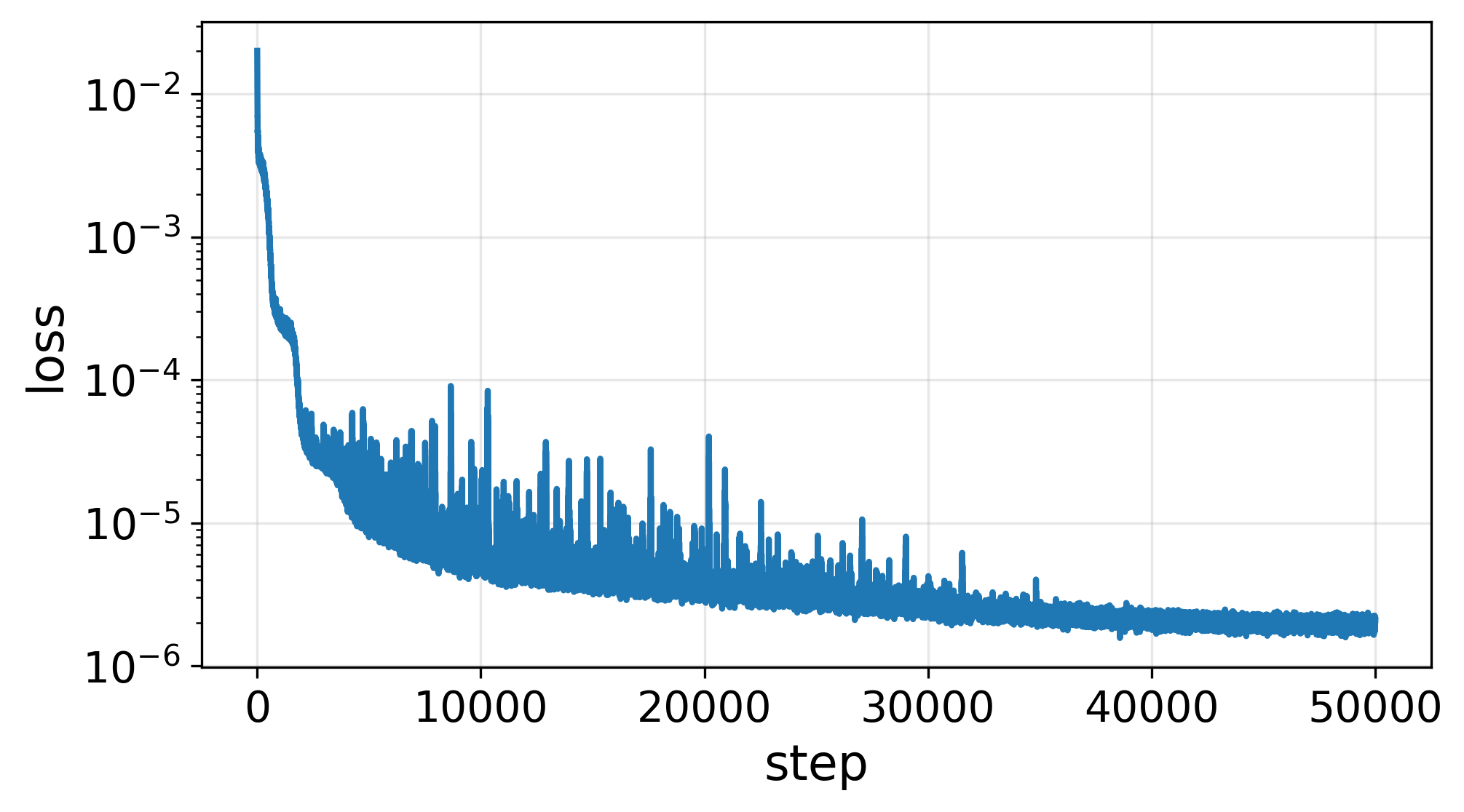}
    \caption{Training loss.}
    \label{fig:test2-loss}
  \end{subfigure}
  \hfill
  \begin{subfigure}{0.48\textwidth}
    \centering
    \includegraphics[width=\textwidth]{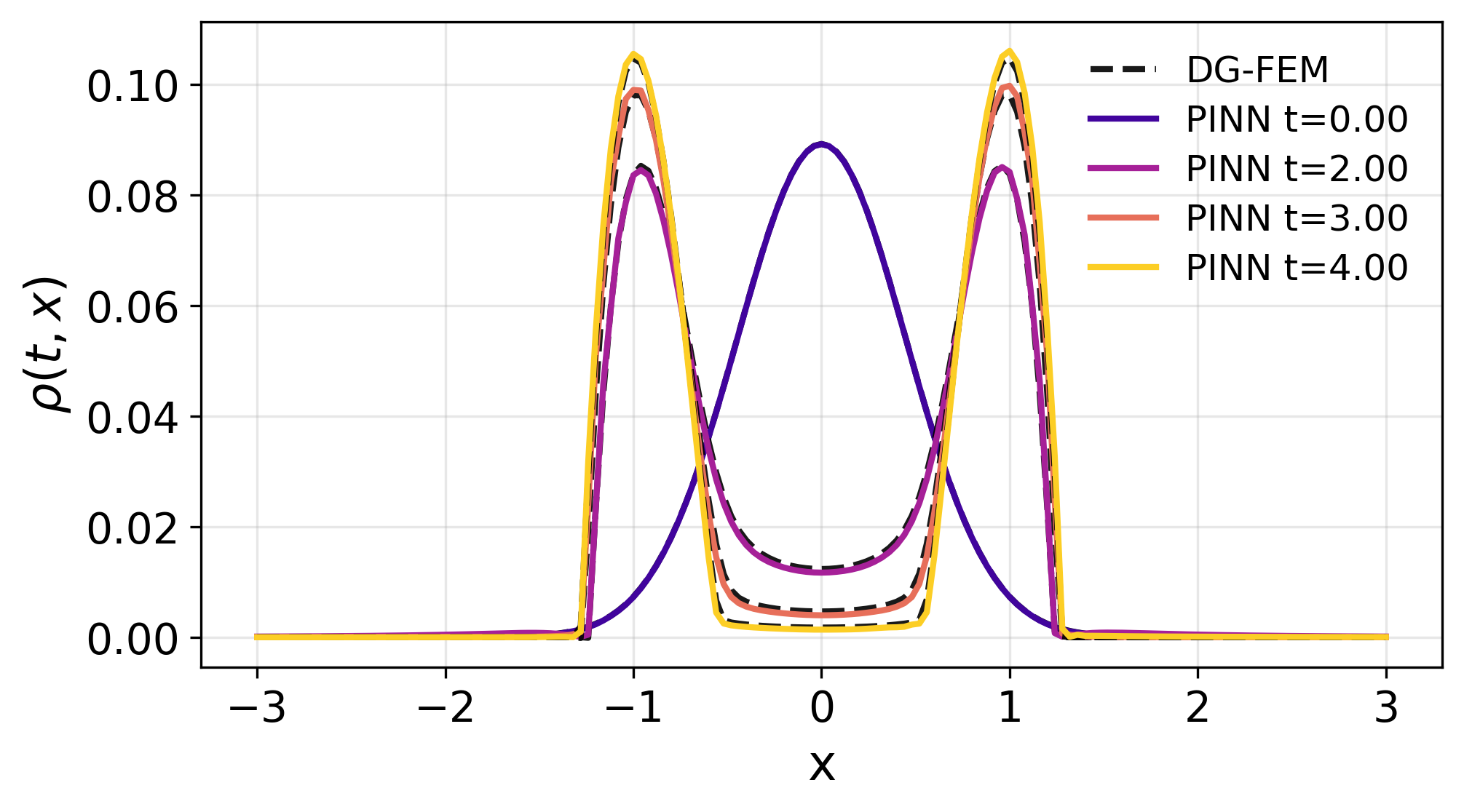}
    \caption{Solution snapshots.}
    \label{fig:test2-snapshots}
  \end{subfigure}

  \vspace{0.4cm}

  \begin{subfigure}{0.48\textwidth}
    \centering
    \includegraphics[width=\textwidth]{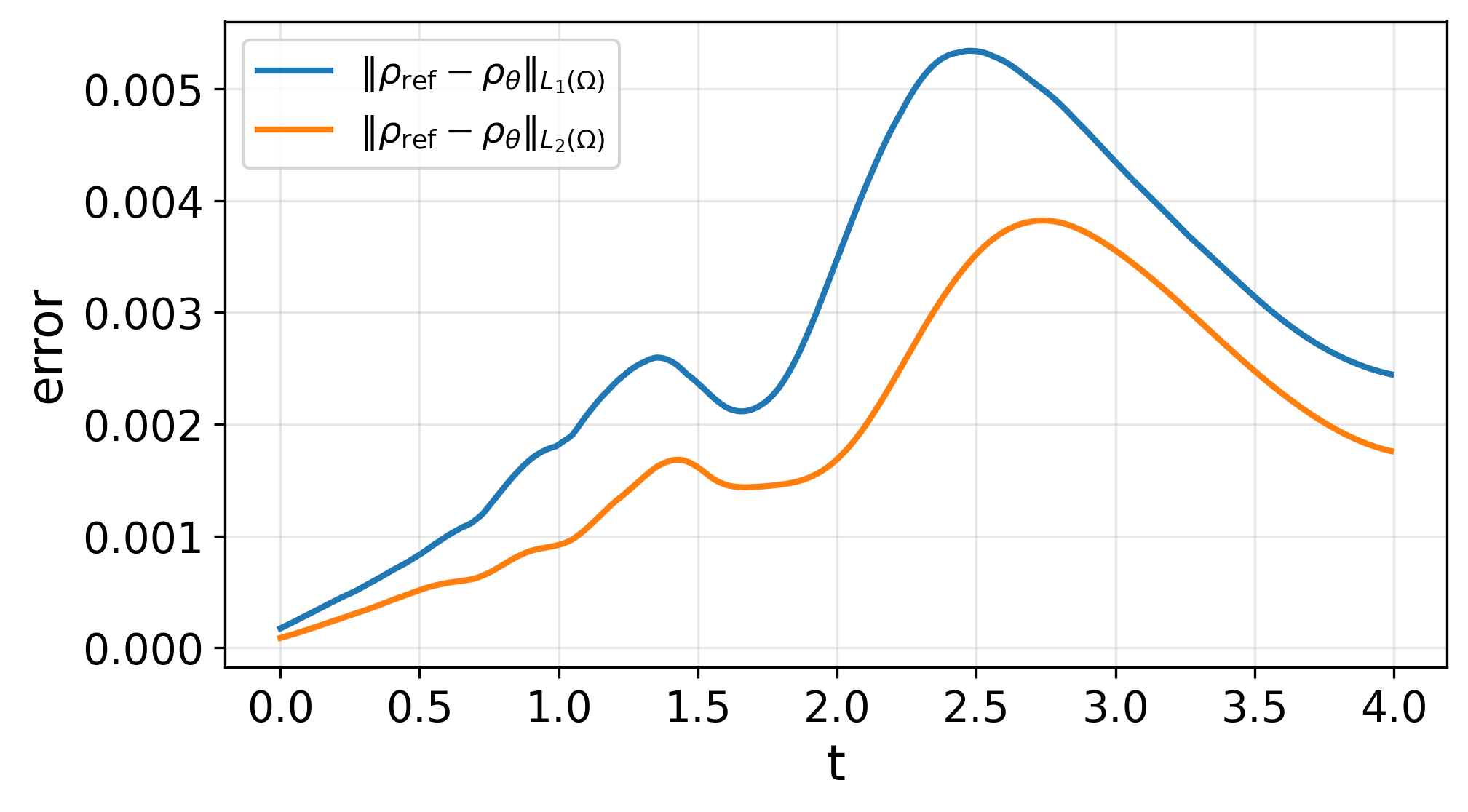}
    \caption{$L^1$ and $L^2$ errors.}
    \label{fig:test2-errors}
  \end{subfigure}
  \hfill
  \begin{subfigure}{0.48\textwidth}
    \centering
    \includegraphics[width=\textwidth]{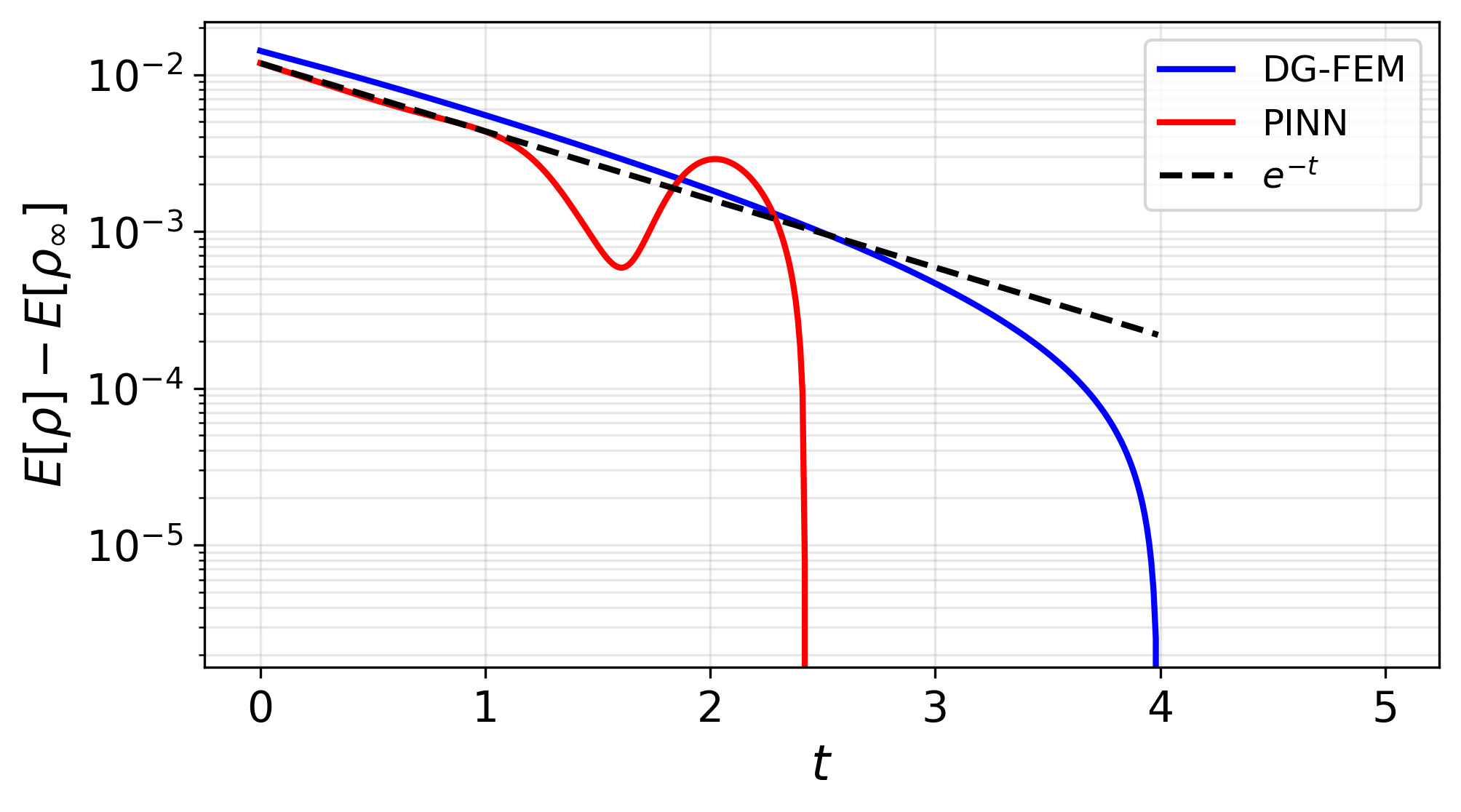}
    \caption{Energy decay.}
    \label{fig:test2-energy}
  \end{subfigure}

  \vspace{0.4cm}

  \begin{subfigure}{0.48\textwidth}
    \centering
    \includegraphics[width=\textwidth]{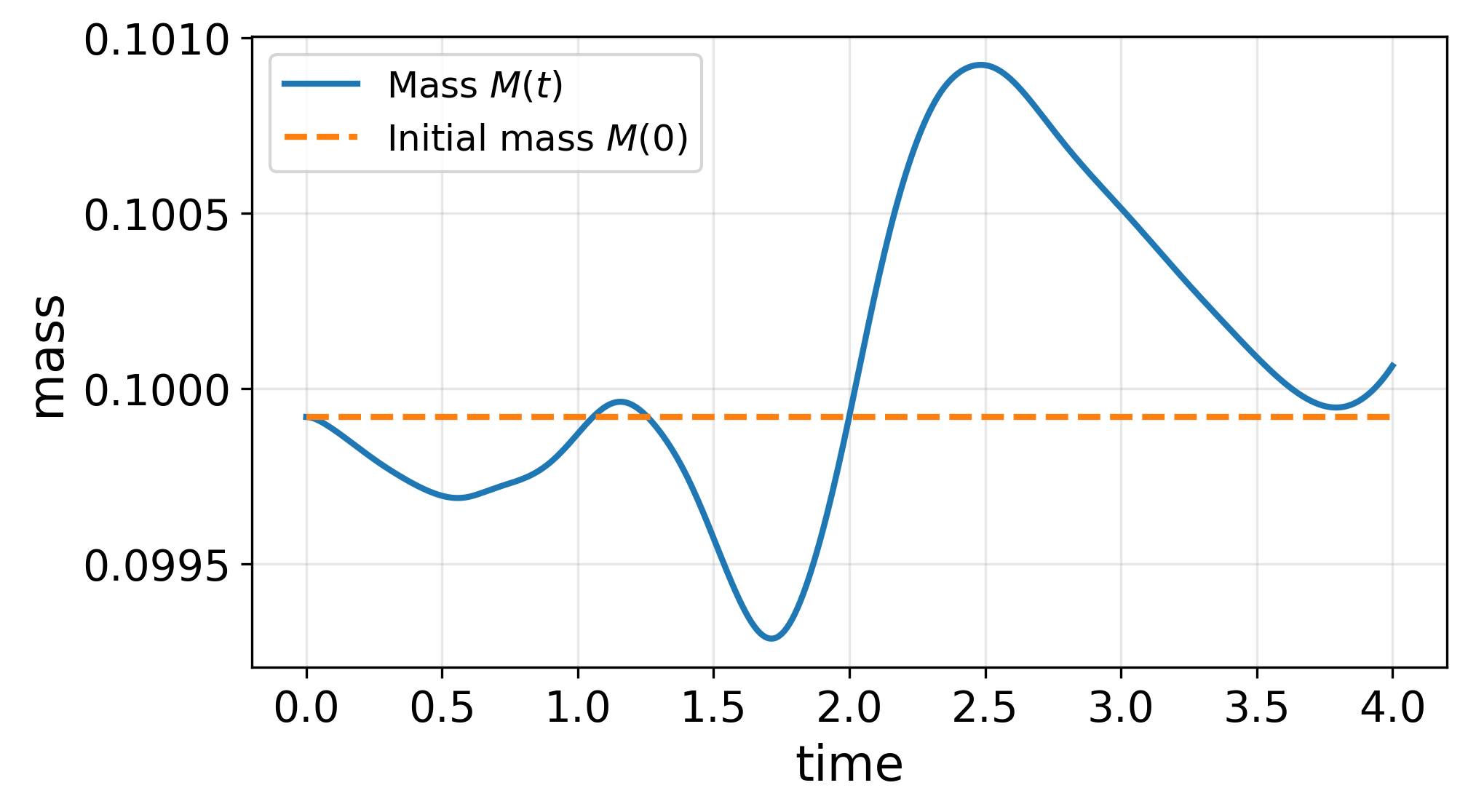}
    \caption{Mass evolution.}
    \label{fig:test2-mass}
  \end{subfigure}
  \hfill
  \begin{subfigure}{0.48\textwidth}
    \centering
    \includegraphics[width=\textwidth]{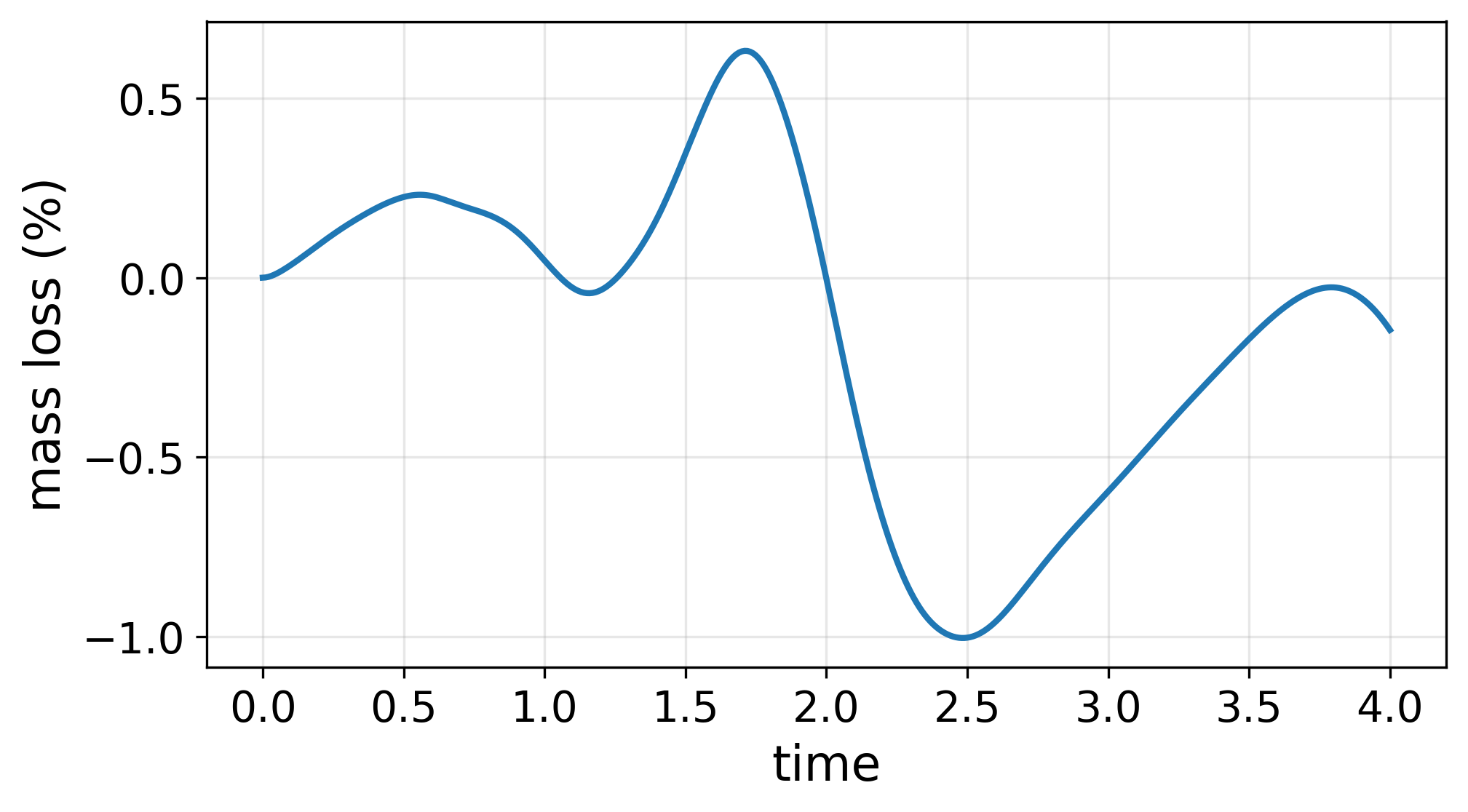}
    \caption{Relative mass loss.}
    \label{fig:test2-massloss}
  \end{subfigure}

  \caption{PINN results for example 2. Top left: training loss. Top right: comparison of the PINN and DG solutions at selected times. Middle left: time evolution of the \(L^1\) and \(L^2\) errors between the PINN and DG-FEM solutions. Middle right: decay of the relative energy \(\E[\rho]-\E[\rho_\infty]\). Bottom left: evolution of the total mass. Bottom right: relative mass loss over time.}
  \label{fig:ex2_pinn}
\end{figure}

Figure~\ref{fig:ex2_pinn} summarizes the numerical results for this example. Panel~(a) shows the decay of the training loss during the trainnig , while panel~(b) compares the neural-network approximation and the reference solution at selected times. The initially centered profile splits into a symmetric bimodal density as time approaches equilibrium, and the agreement between the two solutions remains good throughout the simulated interval. Panel~(c) reports the time evolution of the $L^1$ and $L^2$ errors, which remain small over the full time horizon, although they increase during the transient regime before decreasing slightly near the final time. Panel~(d) shows the decay of the free-energy gap $\E[\rho]-\E[\rho_\infty]$ for both approximations. In both cases, the free energy decreases over time, as expected from the free-energy dissipation property of the problem, although a modest difference in the decay rates becomes visible at late times. Finally, panels~(e) and~(f) report the evolution of the total mass and the relative mass deviation $(M(t)-M_0)/M_0$, respectively, showing that mass is preserved only approximately, with small fluctuations around the initial value over the simulated time interval.

\subsection{Example 3: nonlinear diffusion with nonlocal attractive interactions}
We next consider a one-dimensional nonlinear diffusion problem with nonlocal attractive interaction, given by
\begin{equation*}
\partial_t \rho
=
\partial_x \!\left(
\rho\,\partial_x\left(\nu \rho^{m-1} + (W * \rho)(x)\right)
\right),
\qquad
W(x)=-\max(0,1-|x|),
\end{equation*}
where the interaction kernel $W$ induces attraction through the convolution term $(W * \rho)(x)$.

In this experiment, we set $m=3$ and $\nu=1$, and perform the computations on the domain $\Omega=[-3,3]$ up to the final time $T=4$. The initial condition is  the characteristic function on $[-1.5,\,1.5]$
\begin{equation*}
\rho_0(x)=\dfrac14\,\chi_{[-1.5,\,1.5]}(x).
\end{equation*}
The corresponding free energy is
$$
\E[\rho]=\int_{\Omega}\left(\frac{1}{3}\rho^3+\frac{1}{2}\rho(W * \rho)\right)\,dx.
$$

This test illustrates the interplay between nonlinear diffusion and nonlocal attraction. Starting from a compactly supported piecewise-constant initial density, the solution gradually reorganizes into a single smooth density distribution centered at the origin under the combined influence of the two mechanisms. Similar aggregation--diffusion test cases were considered in \cite{MendesRussoPerezKalliadasis2021}.

For this problem, a single training run over the full time interval did not reliably reproduce the reference solution. To improve stability and accuracy, we adopt a time-windowing strategy in which the interval $[0,T]$ is divided into four consecutive subintervals of equal length. The network is trained sequentially on these subintervals, with the solution obtained at the end of each window used as the initial condition for the next. Using this strategy, the method reproduces the reference evolution accurately, with $1000$ optimization steps performed in each window. We note, however, that time windowing did not improve the results in the previous examples, for which training over the full time interval produced more accurate approximations.

For this experiment, the loss weights are set to
$\lambda_1=1$, $\lambda_2=5$, $\lambda_3=10^{-3}$, and $\lambda_4=0.1$.
The neural network consists of six hidden layers of width $256$ and uses the
modified Lorentzian activation function with $\varepsilon=0.05$. The density
approximation is obtained from the raw network output $g_\theta(x,t)$ using the
identity map. Training is performed with the Adam optimizer, a learning rate of
$10^{-3}$, and a time-batch size of $N_{T,\mathrm{batch}}=20$.

The neural-network approximation is then compared with the reference solution over the full time horizon. The resulting trajectory errors remain moderate and are given by
\begin{align*}
||\rho_{\mathrm{ref}}(x,t)-\rho_{\theta}(x,t)||_{L^2(0,T;L^2(\Omega))}&=  3.04\times 10^{-2},\\
||\rho_{\mathrm{ref}}(x,t)-\rho_{\theta}(x,t)||_{L^1(0,T;L^1(\Omega))}&= 6.83\times 10^{-2}.  
\end{align*}

\begin{figure}[h]
  \centering

  \begin{subfigure}{0.48\textwidth}
    \centering
    \includegraphics[width=\textwidth]{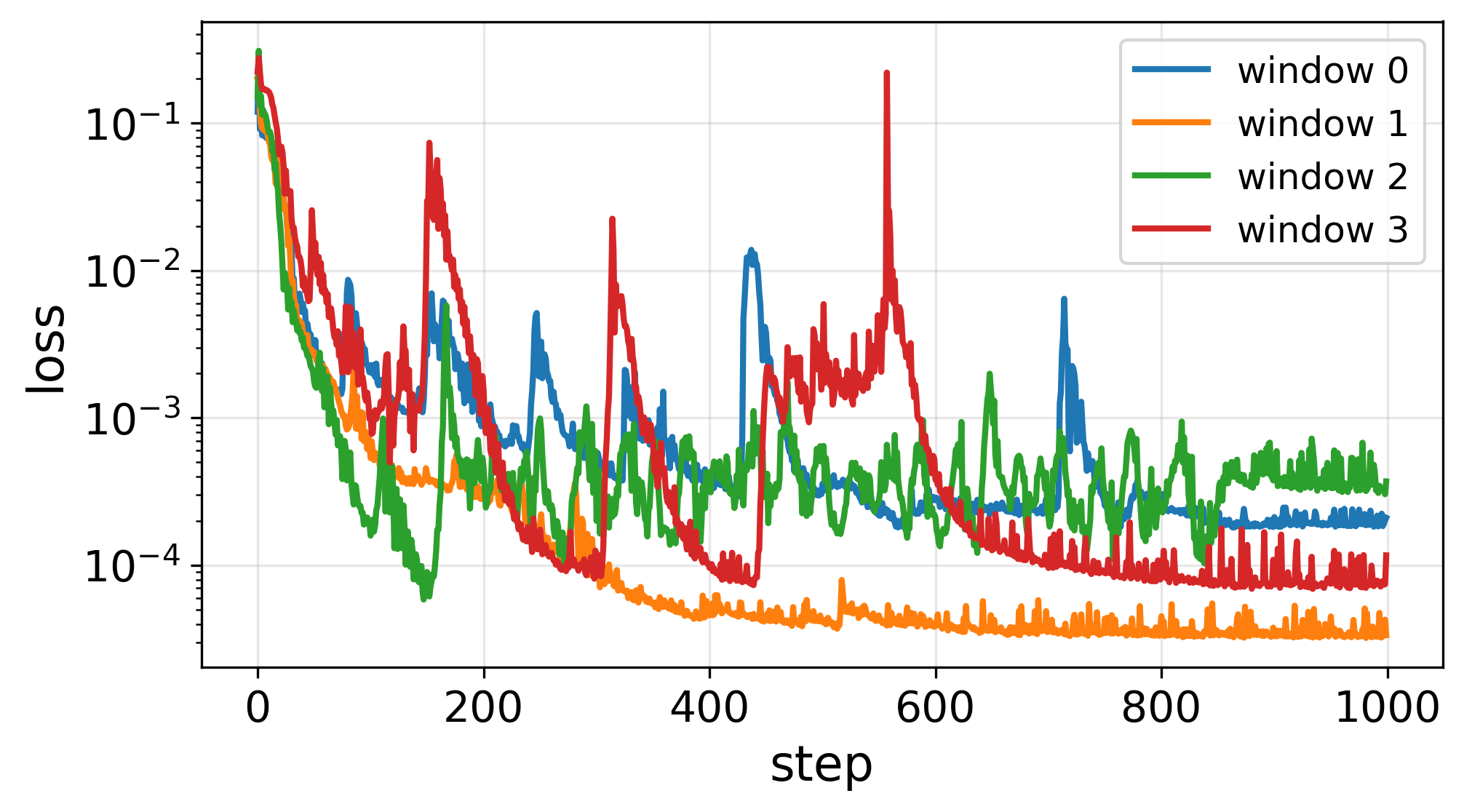}
    \caption{Training loss.}
    \label{fig:test3-loss}
  \end{subfigure}
  \hfill
  \begin{subfigure}{0.48\textwidth}
    \centering
    \includegraphics[width=\textwidth]{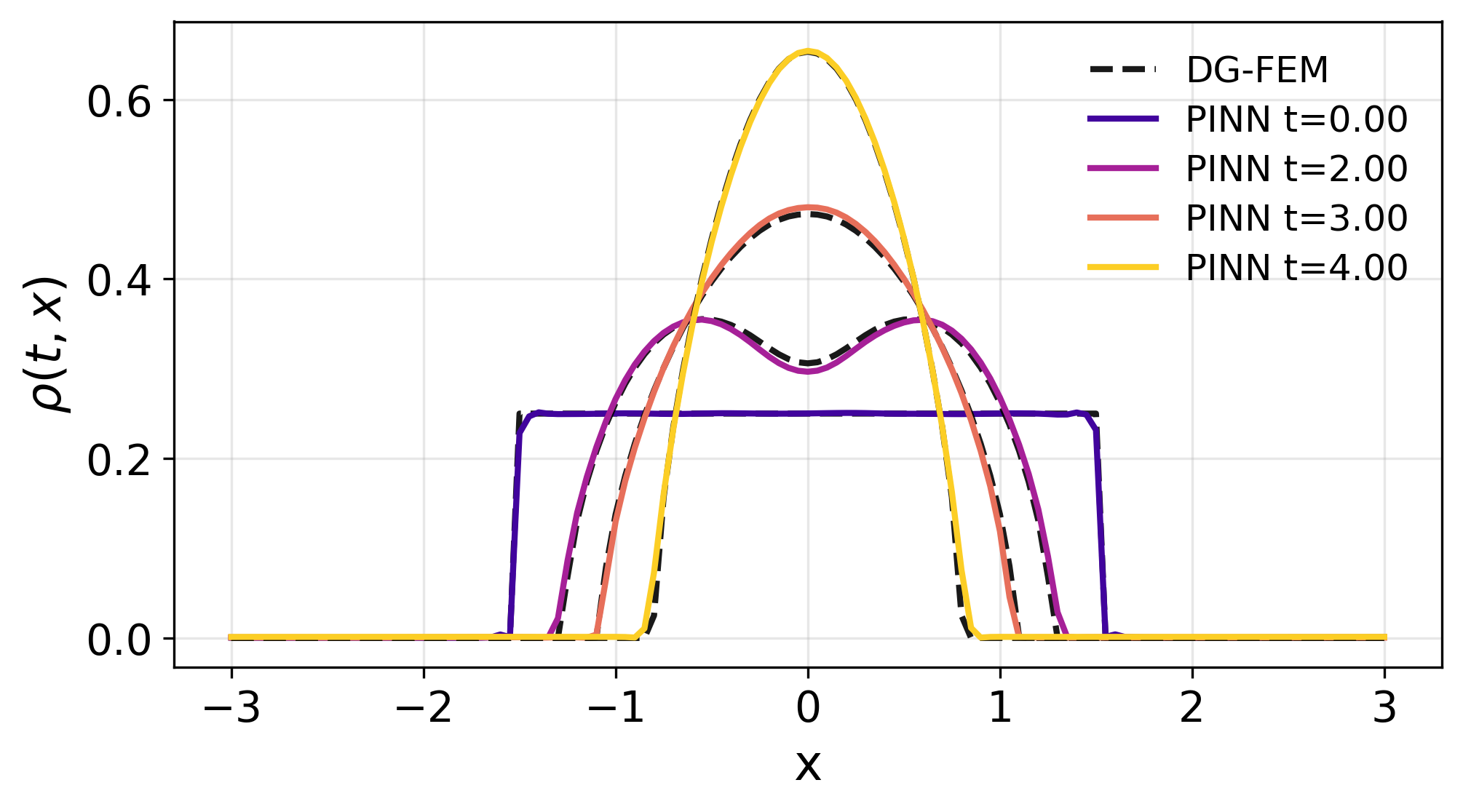}
    \caption{Solution snapshots.}
    \label{fig:test3-snapshots}
  \end{subfigure}

  \vspace{0.4cm}

  \begin{subfigure}{0.48\textwidth}
    \centering
    \includegraphics[width=\textwidth]{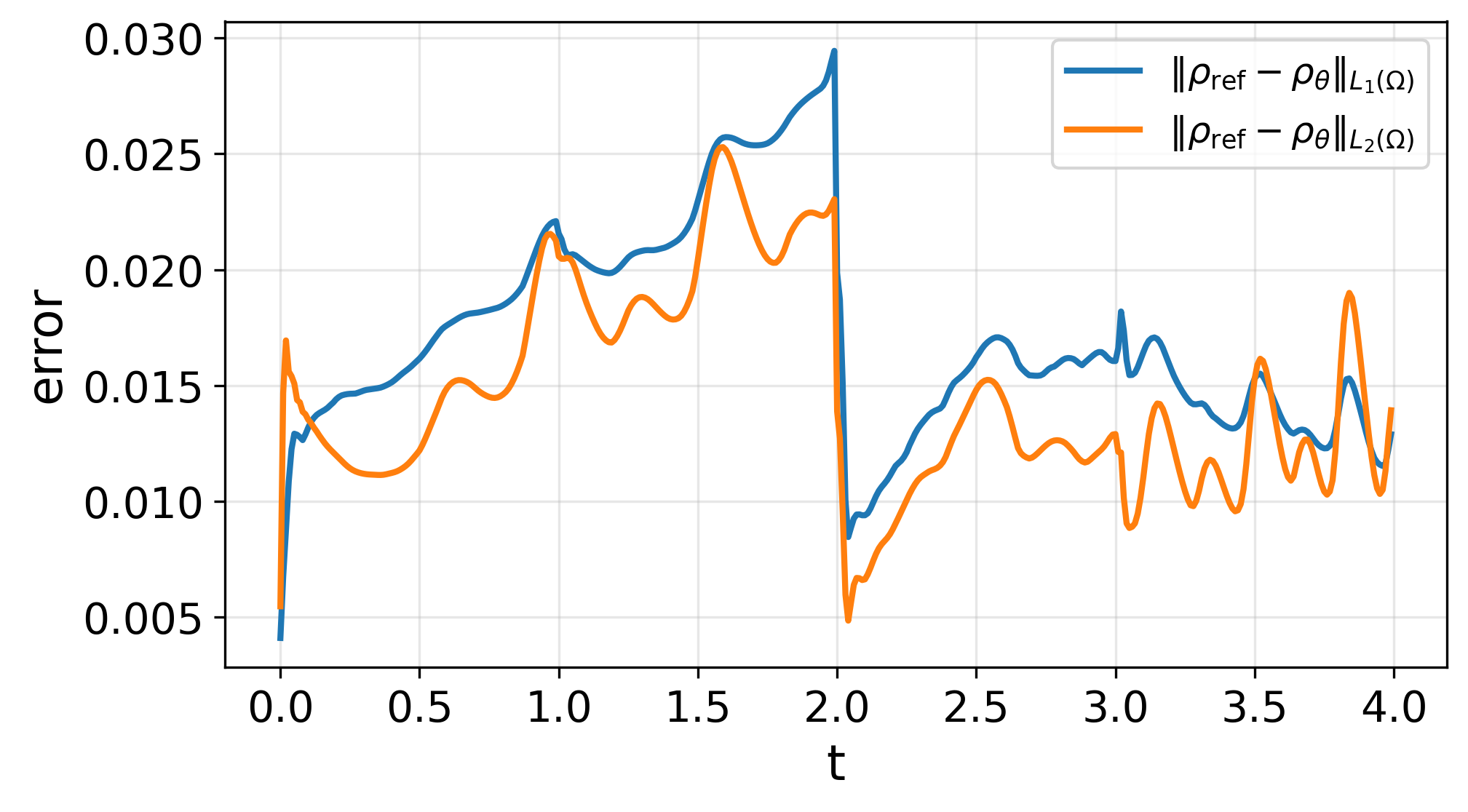}
    \caption{$L^1$ and $L^2$ errors.}
    \label{fig:test3-errors}
  \end{subfigure}
  \hfill
  \begin{subfigure}{0.48\textwidth}
    \centering
    \includegraphics[width=\textwidth]{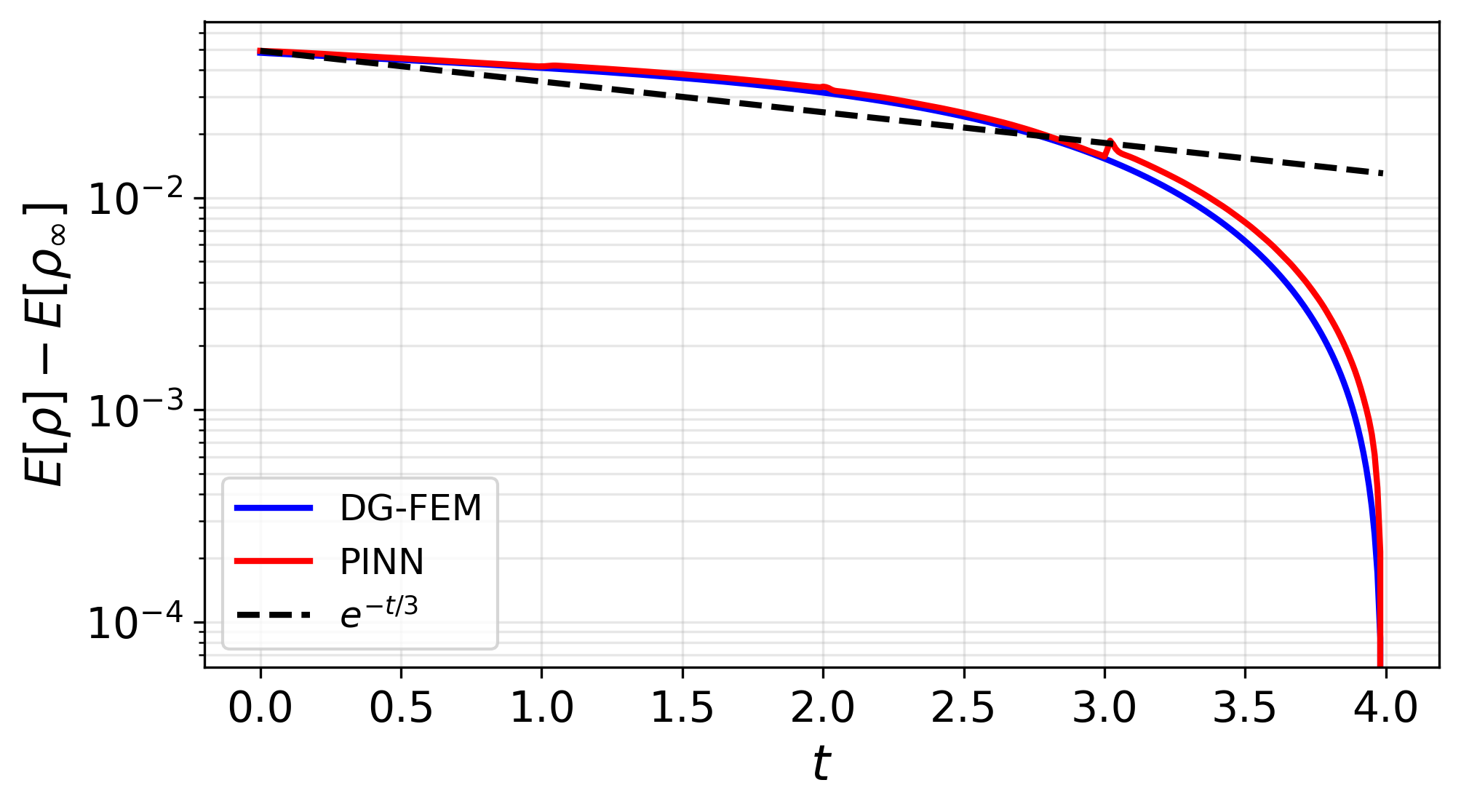}
    \caption{Energy decay.}
    \label{fig:test3-energy}
  \end{subfigure}

  \vspace{0.4cm}

  \begin{subfigure}{0.48\textwidth}
    \centering
    \includegraphics[width=\textwidth]{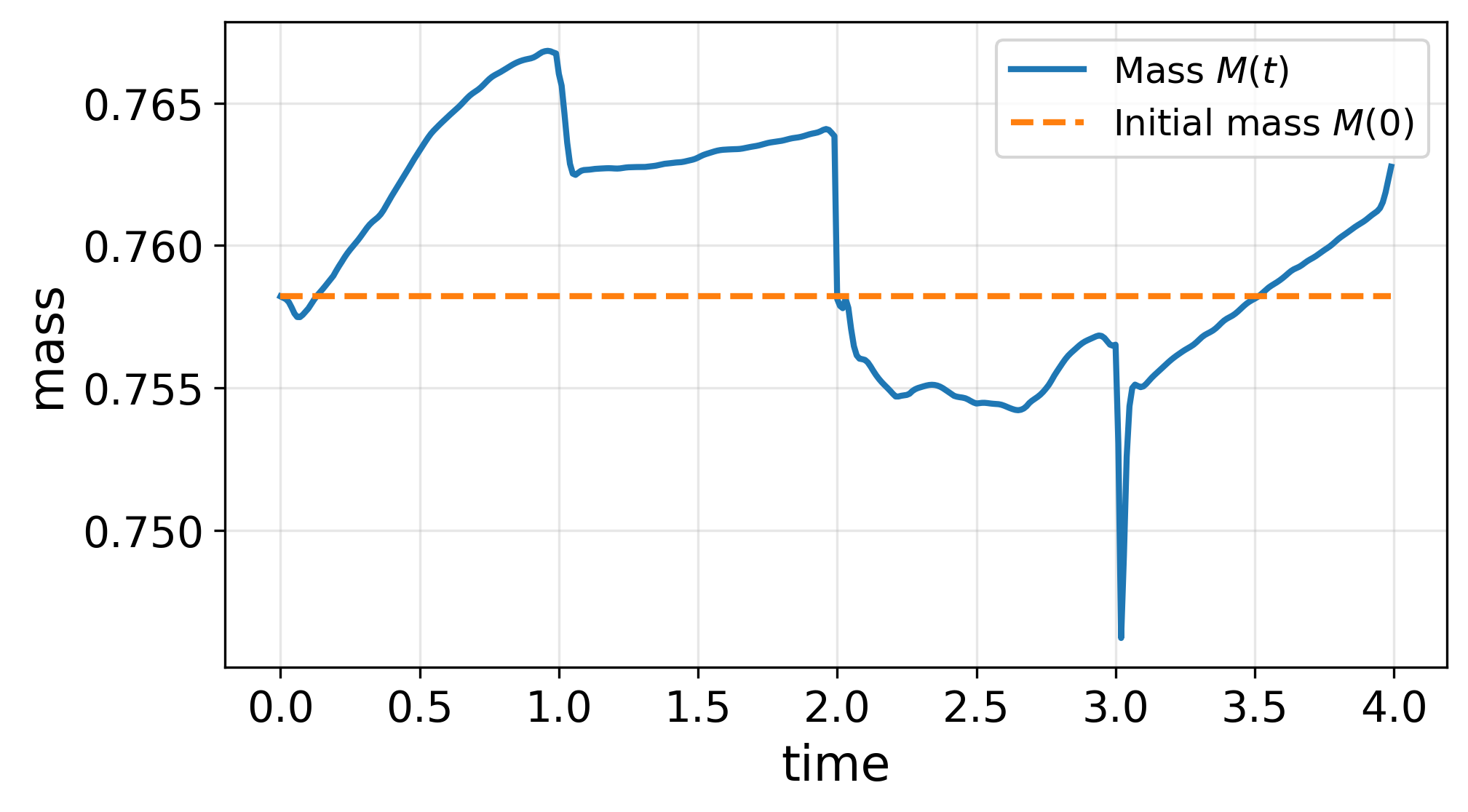}
    \caption{Mass evolution.}
    \label{fig:test3-mass}
  \end{subfigure}
  \hfill
  \begin{subfigure}{0.48\textwidth}
    \centering
    \includegraphics[width=\textwidth]{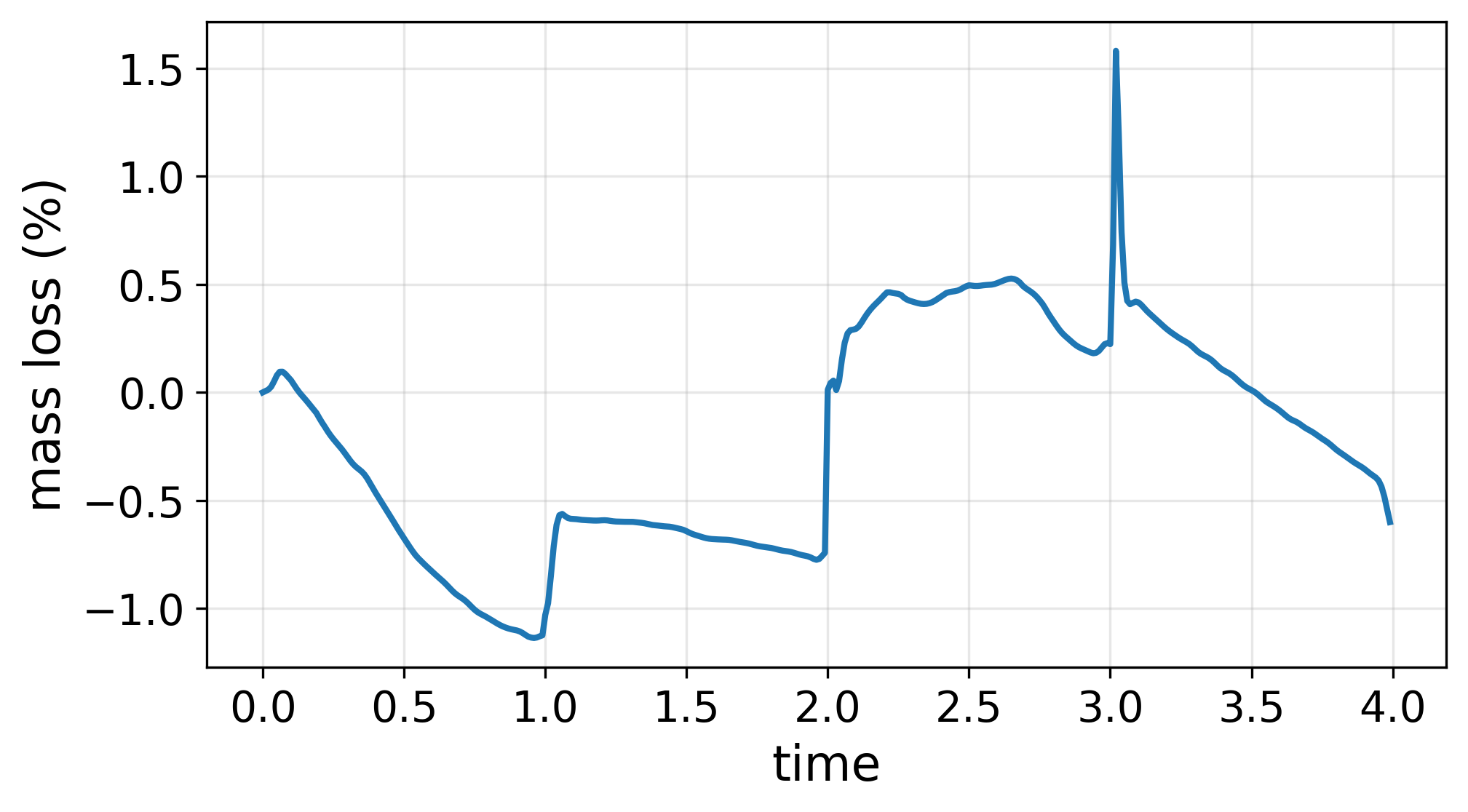}
    \caption{Relative mass loss.}
    \label{fig:test3-massloss}
  \end{subfigure}

\caption{PINN results for example 3 using a sequential time-window training strategy. Top left: training loss for each time window. Top right: comparison of the neural-network approximation and the DG reference solution at selected times. Middle left: time evolution of the \(L^1\) and \(L^2\) errors between the two solutions. Middle right: decay of the free-energy gap \(\E[\rho]-\E[\rho_\infty]\). Bottom left: evolution of the total mass. Bottom right: relative mass deviation over time.}
\label{fig:ex3_pinn}
\end{figure}

Figure~\ref{fig:ex3_pinn} summarizes the numerical results. Panel~(a) shows the training loss in each time window. Panel~(b) compares the neural-network and reference solutions at selected times and shows the evolution of the initially compactly supported density into a single smooth distribution centered at the origin. Panel~(c) reports the $L^1$ and $L^2$ errors over time. Changes at the window interfaces arise from the sequential training procedure. Panel~(d) shows the free-energy gap $\E[\rho]-\E[\rho_\infty]$ for both solutions. In each case, the free energy decreases monotonically, as expected from the gradient-flow formulation, and the two curves remain close throughout the simulation. Panels~(e) and~(f) show the total mass and the relative mass deviation, respectively. The changes at the window interfaces are caused by the transition between consecutive training intervals. Overall, the method reproduces the reference evolution and energy decay well, although the time-windowing procedure introduces small discontinuities in the error and mass curves.

\subsection{Example 4: droplet coalescence without external confinement}

We next consider a two-dimensional model of droplet coalescence without external confinement. The free energy combines ideal-gas entropy, local short-range repulsion, and nonlocal attractive interactions. The model follows the classical density-functional formulation considered, for example, in \cite{YatsyshinKalliadasis2021}.

The density evolves according to
\begin{equation}
\partial_t \rho
=
\nabla \cdot \left(
\rho \nabla
\left(
H(\rho)
+
V(\mathbf{x})
+
(W * \rho)(\mathbf{x})
\right)
\right),
\qquad
\text{in }\Omega\times(0,T],
\label{eq:ddft_coalescence_pde}
\end{equation}
where the chemical potential is
\begin{equation}
\mu(\mathbf{x},t)
=
k_B\tau\ln\rho(\mathbf{x},t)
+
\psi\bigl(\rho(\mathbf{x},t)\bigr)
+
\rho(\mathbf{x},t)\psi'\bigl(\rho(\mathbf{x},t)\bigr)
+
\bigl(\varphi_{\mathrm{attr}}*\rho\bigr)(\mathbf{x},t).
\label{eq:ddft_coalescence_mu}
\end{equation}
The problem is supplemented by the no-flux boundary condition~\eqref{eq:noflux}.

In the notation of \eqref{eq:energy}, this corresponds to
\begin{equation}
H(\rho)
=
k_B\tau\,\rho(\ln\rho-1)
+
\rho\,\psi(\rho),
\qquad
V(\mathbf{x})=0,
\qquad
W(\mathbf{x}-\mathbf{y})
=
\varphi_{\mathrm{attr}}(|\mathbf{x}-\mathbf{y}|).
\label{eq:coalescence_HVW}
\end{equation}

The local repulsive contribution is described by the Carnahan--Starling expression \cite{Lutsko2010},
\begin{equation*}
\psi(\rho)
=
k_B\tau\,
\dfrac{\eta(4-3\eta)}{(1-\eta)^2},
\qquad
\eta
=
\dfrac{\pi\sigma^3\rho}{6},
\end{equation*}
where $k_B$ is the Boltzmann constant, $\tau$ is the temperature, and $\sigma$ denotes the particle diameter.

The nonlocal attractive interaction is defined by
\begin{equation*}
\varphi_{\mathrm{attr}}(r)
=
\begin{cases}
0, & r\leq\sigma,\\[4pt]
\varphi^{12-6}_{\epsilon,\sigma}(r), & r>\sigma,
\end{cases}
\qquad
\varphi^{12-6}_{\epsilon,\sigma}(r)
=
4\epsilon
\left[
\left(\frac{\sigma}{r}\right)^{12}
-
\left(\frac{\sigma}{r}\right)^6
\right],
\end{equation*}
where $\epsilon$ controls the strength of the attractive interaction.

The corresponding free energy is
\begin{equation*}
\begin{aligned}
\E[\rho]
&=
\int_{\Omega}
\left[
k_B\tau\,\rho(\mathbf{x})
\left(\ln\rho(\mathbf{x})-1\right)
+
\rho(\mathbf{x})\,
\psi\left(\rho(\mathbf{x})\right)
\right]
\,d\mathbf{x}
\\
&\quad+
\frac{1}{2}
\int_{\Omega}\int_{\Omega}
\rho(\mathbf{x})\rho(\mathbf{y})
\varphi_{\mathrm{attr}}
\left(|\mathbf{x}-\mathbf{y}|\right)
\,d\mathbf{y}\,d\mathbf{x}.
\end{aligned}
\end{equation*}
The three contributions represent the ideal-gas free energy, local repulsion, and nonlocal attraction, respectively.

The initial density represents the smooth union of two diffuse droplets:
\begin{equation*}
\rho_0(\mathbf{x})
=
\rho_{\text{g}}
+
(\rho_{\text{l}}-\rho_{\text{g}})\,s(\mathbf{x}),
\end{equation*}
where $\rho_{\text{g}}$ and $\rho_{\text{l}}$ denote the gas and liquid densities, respectively, and
\begin{equation*}
s(\mathbf{x})
=
1-
\left(1-s_1(\mathbf{x})\right)
\left(1-s_2(\mathbf{x})\right),
\end{equation*}
with
\begin{equation*}
s_j(\mathbf{x})
=
\frac{1}{2}
\left[
1-
\tanh\left(
\frac{|\mathbf{x}-\mathbf{c}_j|-R_j}{w}
\right)
\right],
\qquad
j=1,2.
\end{equation*}
Here, $\mathbf{c}_j$ and $R_j$ denote the centre and radius of droplet $j$, respectively, while $w>0$ controls the interfacial width. In the computations, we use
\begin{equation*}
\rho_g=10^{-3},
\qquad
\rho_l=1,
\qquad
\mathbf{c}_1=(7,10),
\qquad
\mathbf{c}_2=(13,10),
\qquad
R_1=R_2=2,
\qquad
w=1.
\end{equation*}

The initial state therefore consists of two diffuse droplets surrounded by a low-density gas. The nonlocal attraction causes the droplets to approach one another and merge into a single droplet.

For this experiment, the loss weights are set to
$\lambda_1=1$, $\lambda_2=10$, $\lambda_3=1$, and
$\lambda_4=\textcolor{red}{\text{?}}$.
The computations are performed on the domain
$\Omega=[0,20]\times[0,20]$ up to the final time $T=7$, using a
$40\times40$ Cartesian grid. The temperature is set to $\tau=0.4$.
The neural network uses the modified Lorentzian activation function with
$\varepsilon=0.05$, and the density approximation is obtained from the raw
network output $g_\theta(\mathbf{x},t)$ using the identity map. Training is
performed with the Adam optimizer, a learning rate of $10^{-3}$, a time-batch
size of $100$, and a total of $10{,}000$ optimization steps.

\begin{figure}[htbp]
  \centering

  \begin{subfigure}{0.48\textwidth}
    \centering
    \includegraphics[width=\textwidth]{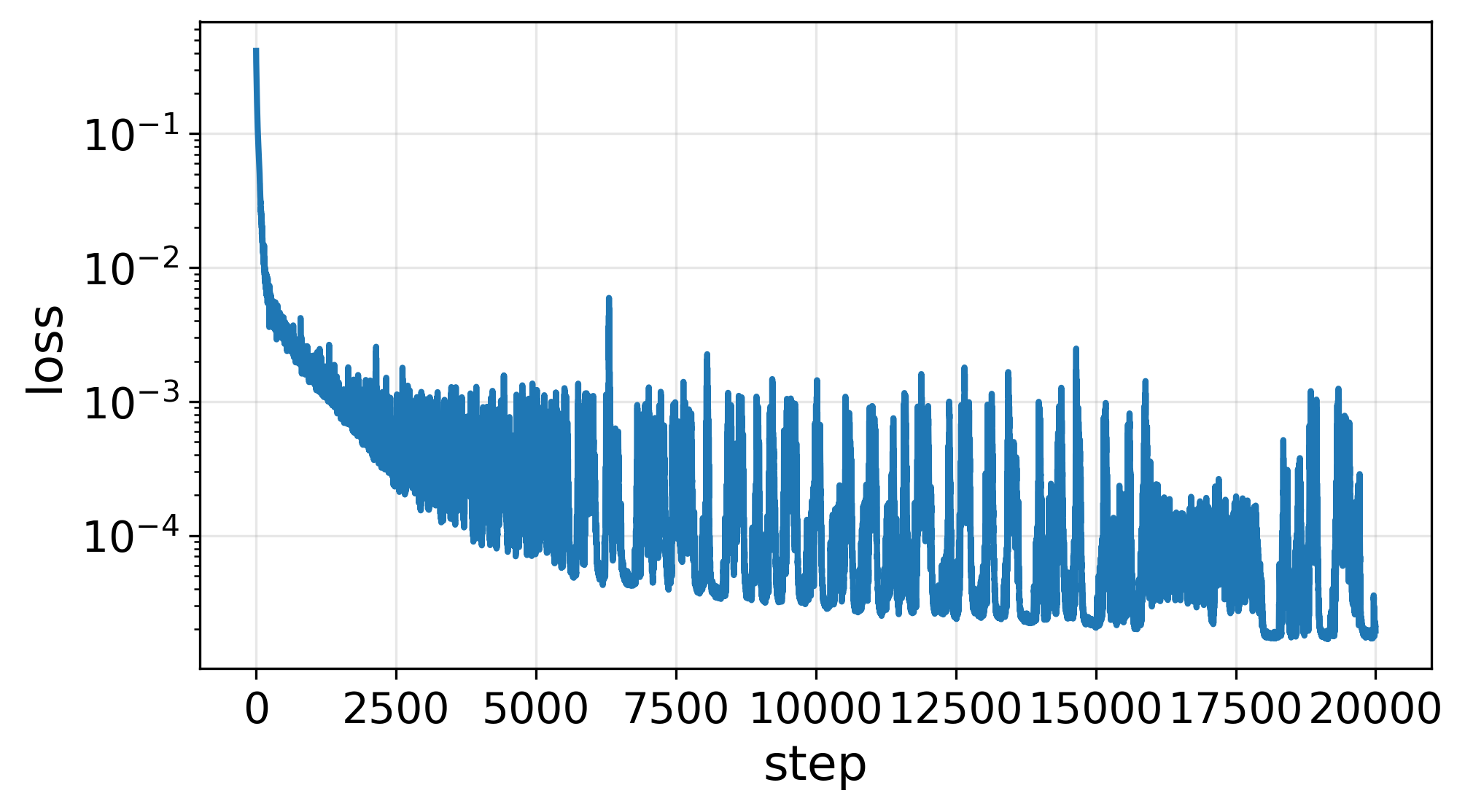}
    \caption{Training loss.}
    \label{fig:2d-loss}
  \end{subfigure}
  \hfill
  \begin{subfigure}{0.48\textwidth}
    \centering
    \includegraphics[width=\textwidth]{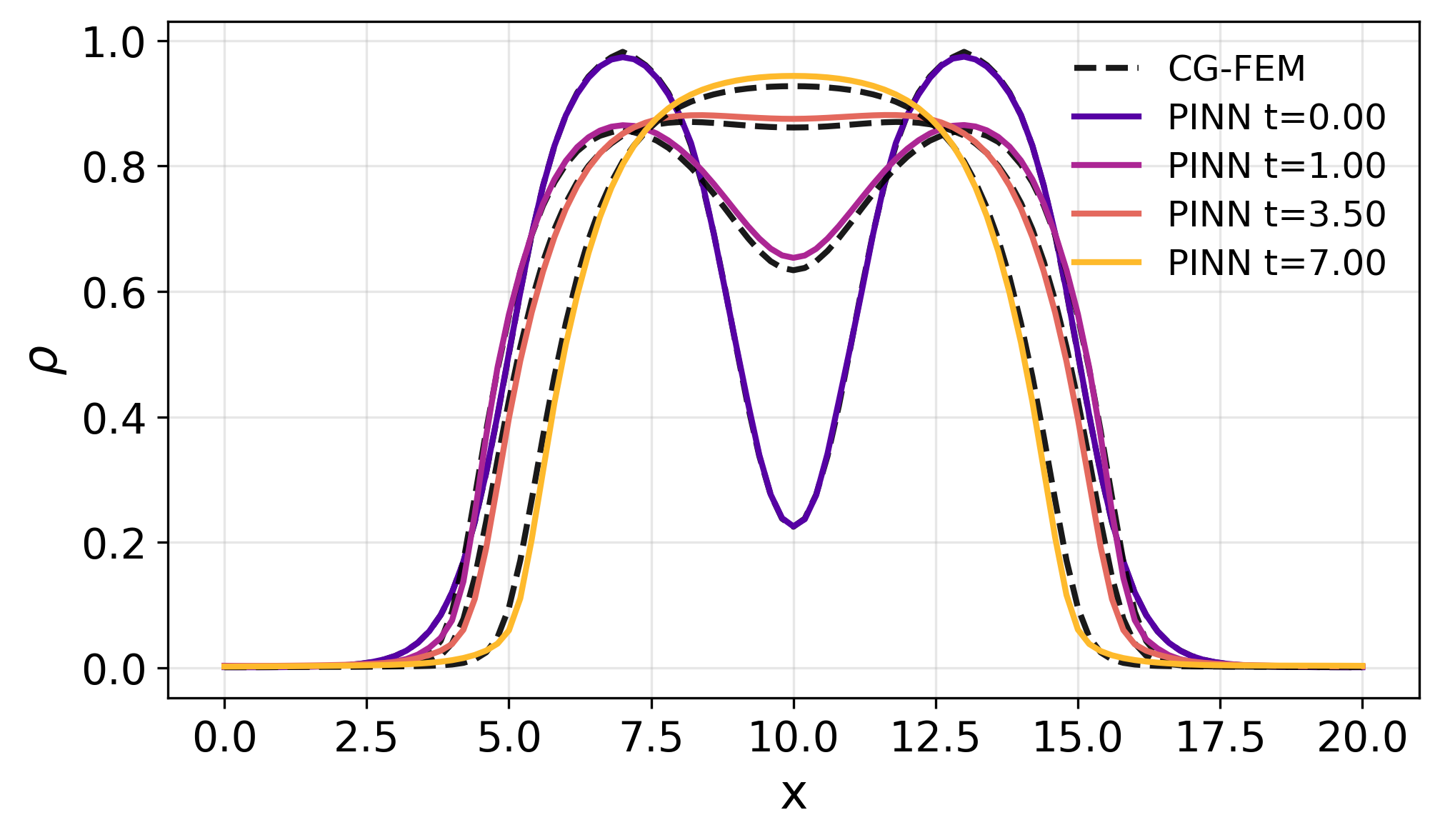}
    \caption{Solution at $y=10$ for selected times.}
    \label{fig:2d-cuts}
  \end{subfigure}

  \vspace{0.4cm}

  \begin{subfigure}{0.48\textwidth}
    \centering
    \includegraphics[width=\textwidth]{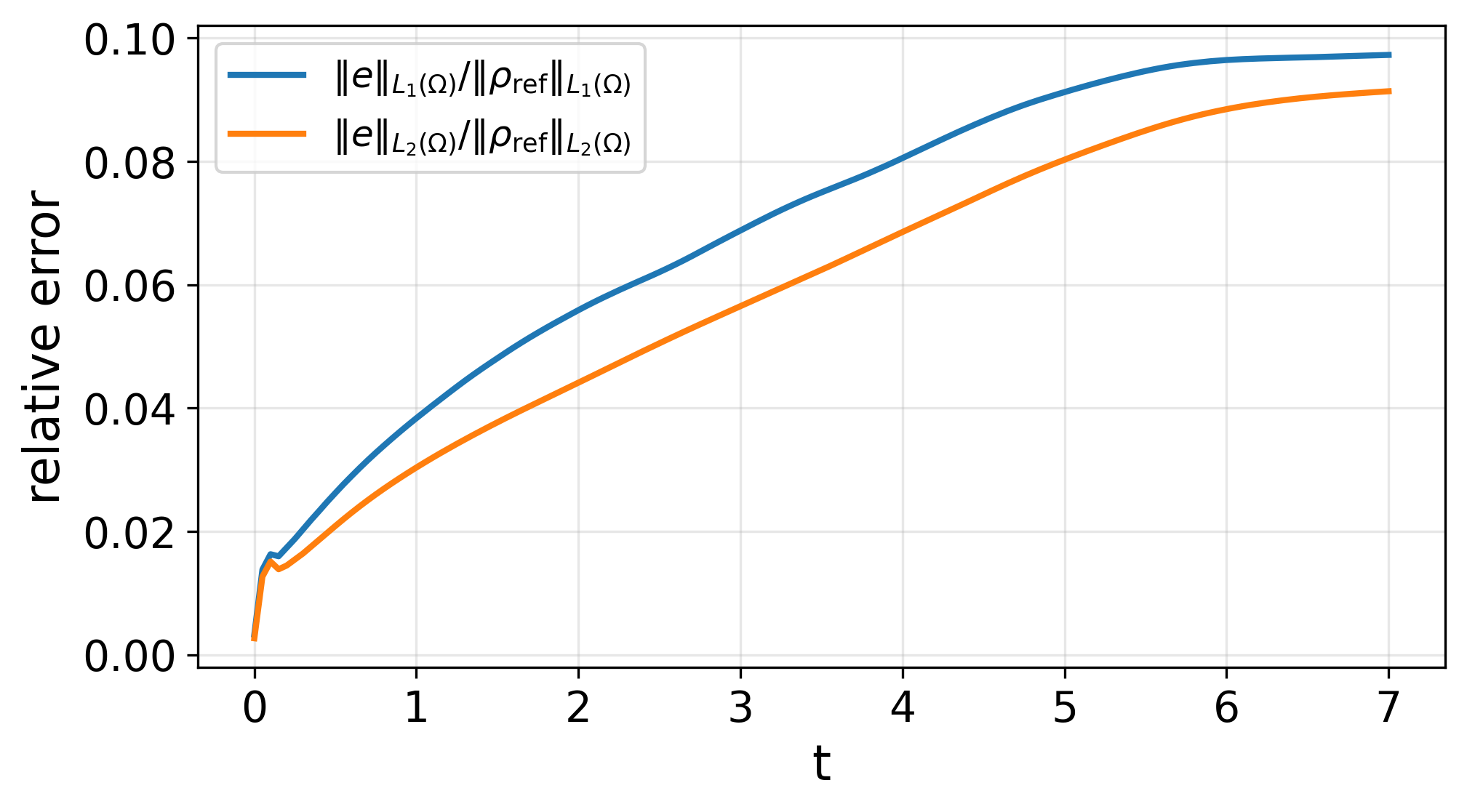}
    \caption{Relative $L^1$ and $L^2$ errors versus time.}
    \label{fig:2d-errors}
  \end{subfigure}
  \hfill
  \begin{subfigure}{0.48\textwidth}
    \centering
    \includegraphics[width=\textwidth]{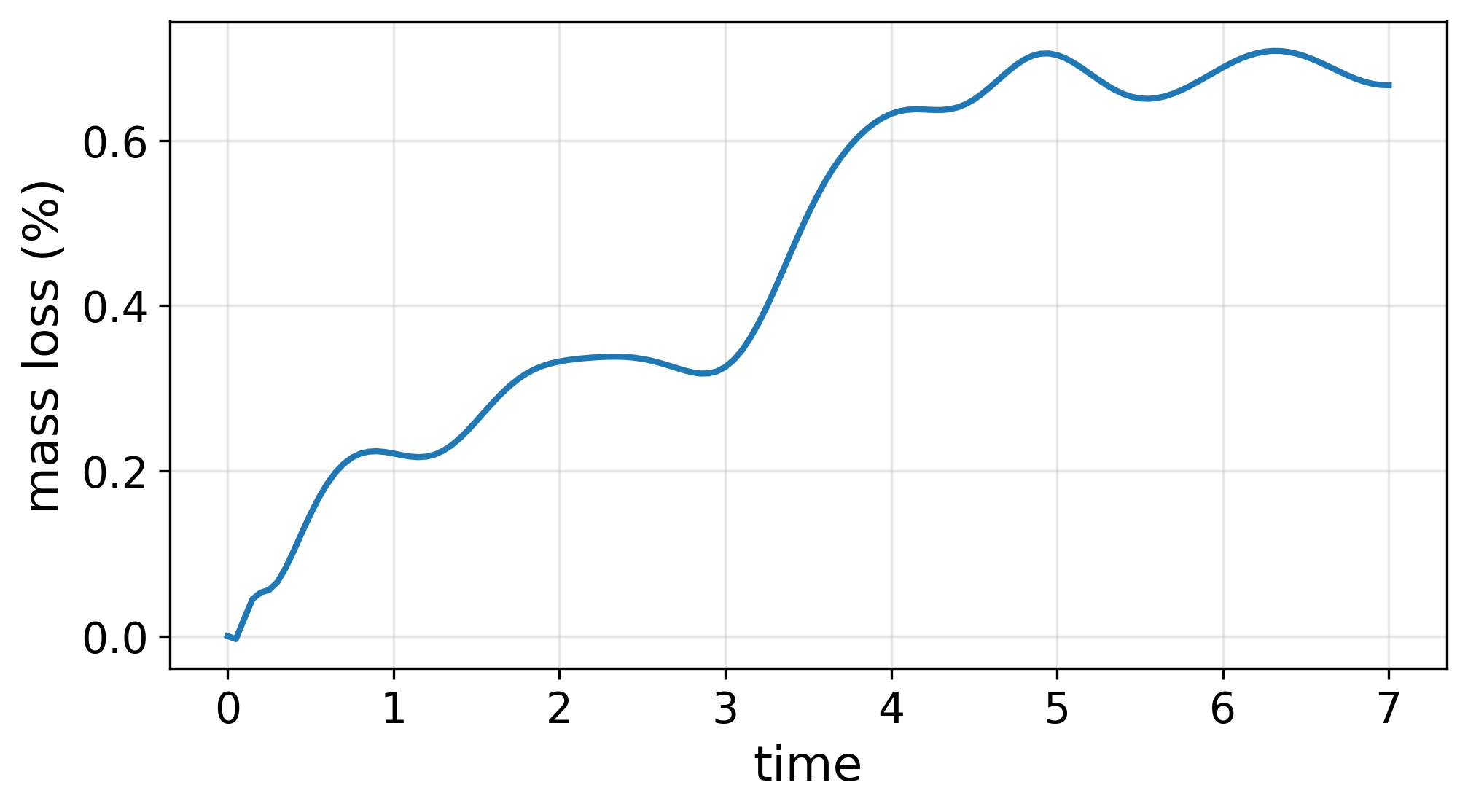}
    \caption{Relative mass loss of the PINN solution.}
    \label{fig:2d-massloss}
  \end{subfigure}

  \caption{Two-dimensional comparison between the PINN and the FEM reference solution. The figure shows the training loss, one-dimensional cross-sections at $y=10$, the relative $L^1$ and $L^2$ errors over time, and the relative mass loss of the PINN solution.}
  \label{fig:test5-2d-summary}
\end{figure}

Figure~\ref{fig:test5-2d-summary} summarizes the overall performance of the PINN for the two-dimensional droplet coalescence problem. Panel~(a) shows a decrease of the training loss, with fluctuations due to the stochastic nature of mini-batch training. Larger batch sizes may reduce these oscillations, but at a higher memory cost. Panel~(b) compares one-dimensional cross-sections along the line $y=10$, viewed as functions of $x$, at selected times. This panel shows that the PINN reproduces the CG-FEM reference solution well, with the largest discrepancy occurring during the transient stage of coalescence. Panel~(c) reports the relative $L^1$ and $L^2$ errors over time. Both errors remain moderate, but increase gradually during the evolution, indicating that small discrepancies accumulate as the solution progresses toward the merged configuration, before stabilizing as the solution approaches equilibrium. Panel~(d) shows the relative mass loss of the PINN solution. Although the mass is not preserved exactly, the deviation remains below $1\%$ over the full time interval.

\begin{figure}[htbp]
  \centering

  \begin{subfigure}{\textwidth}
    \centering
    \includegraphics[width=0.32\textwidth]{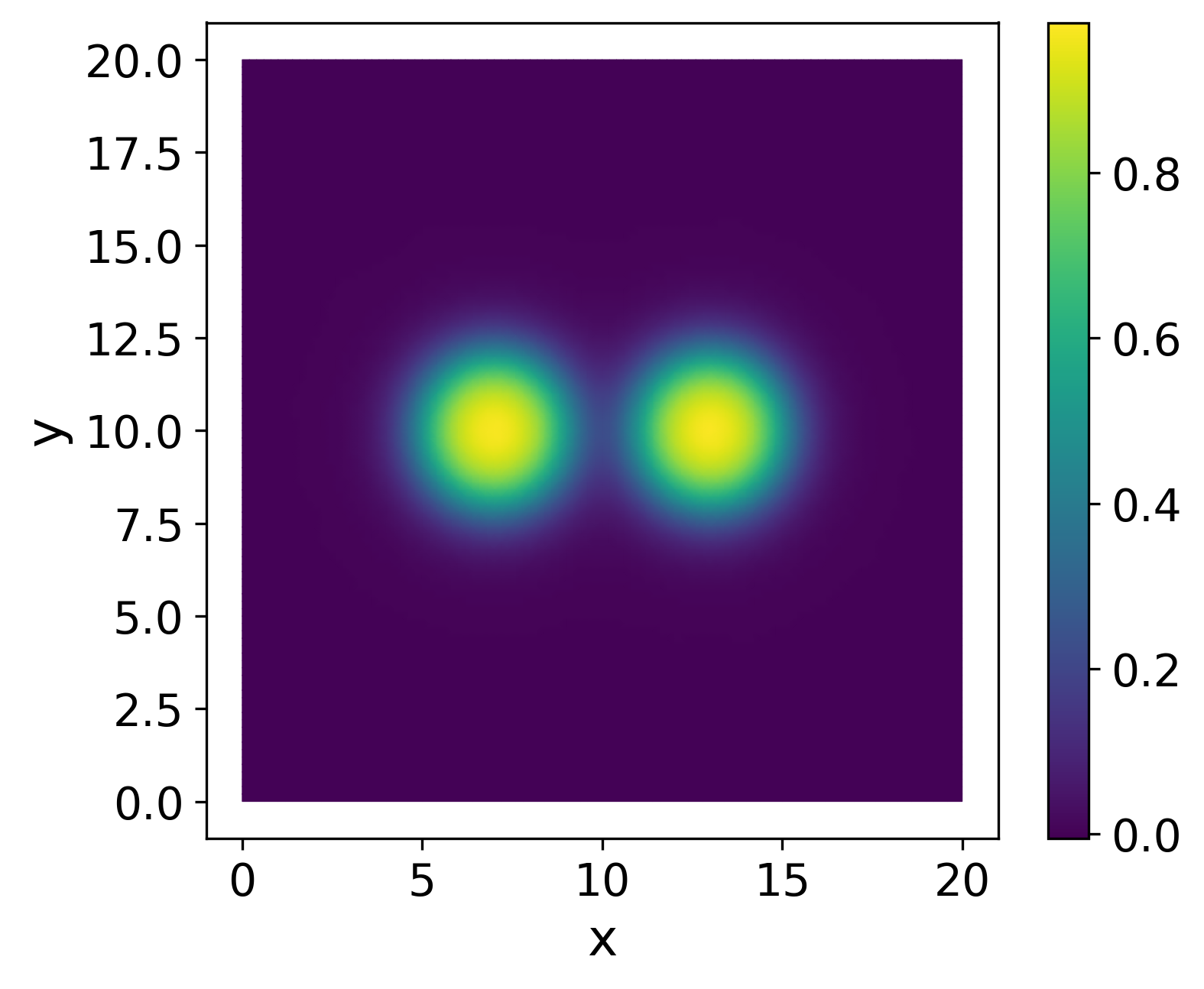}
    \hfill
    \includegraphics[width=0.32\textwidth]{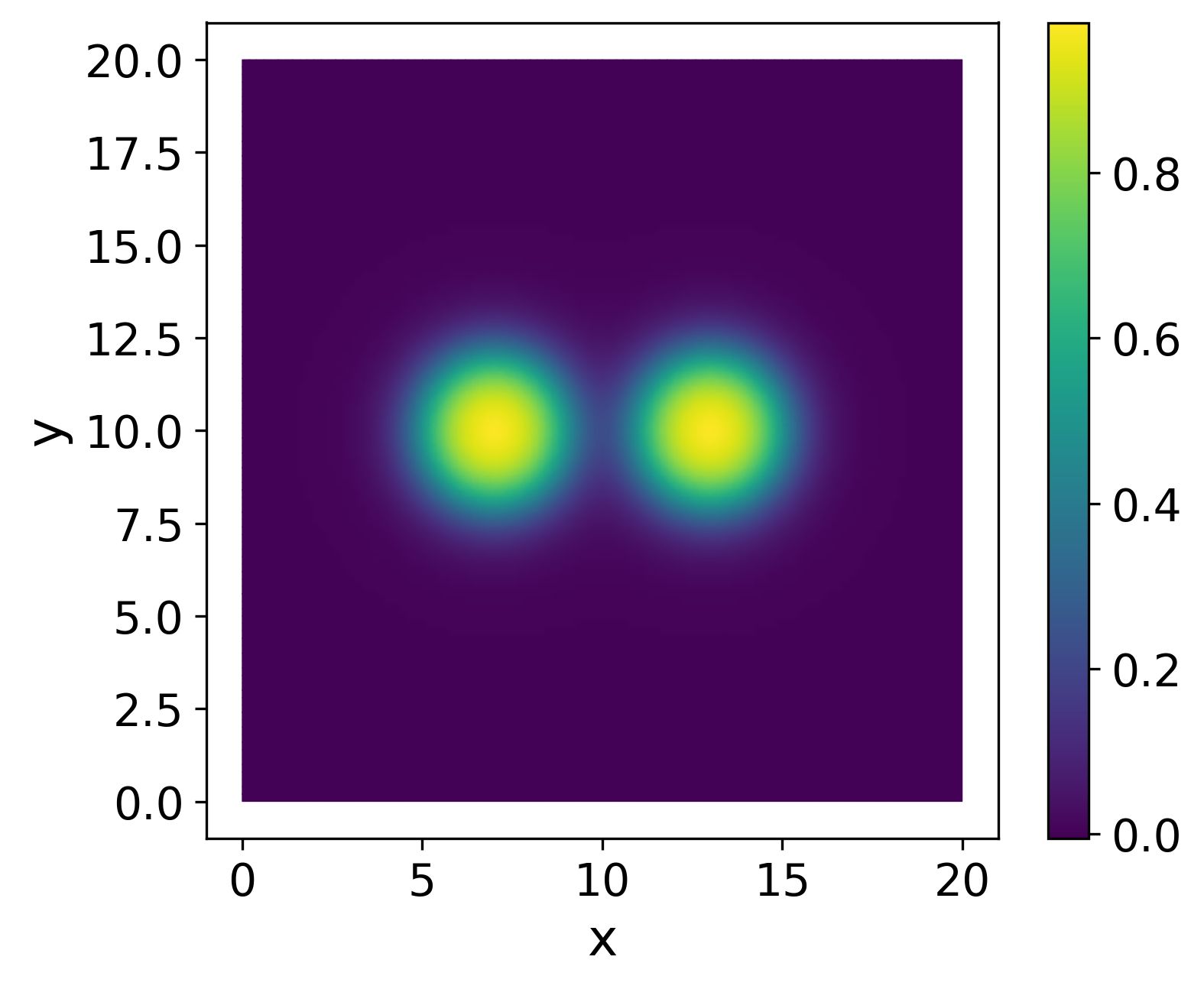}
    \hfill
    \includegraphics[width=0.32\textwidth]{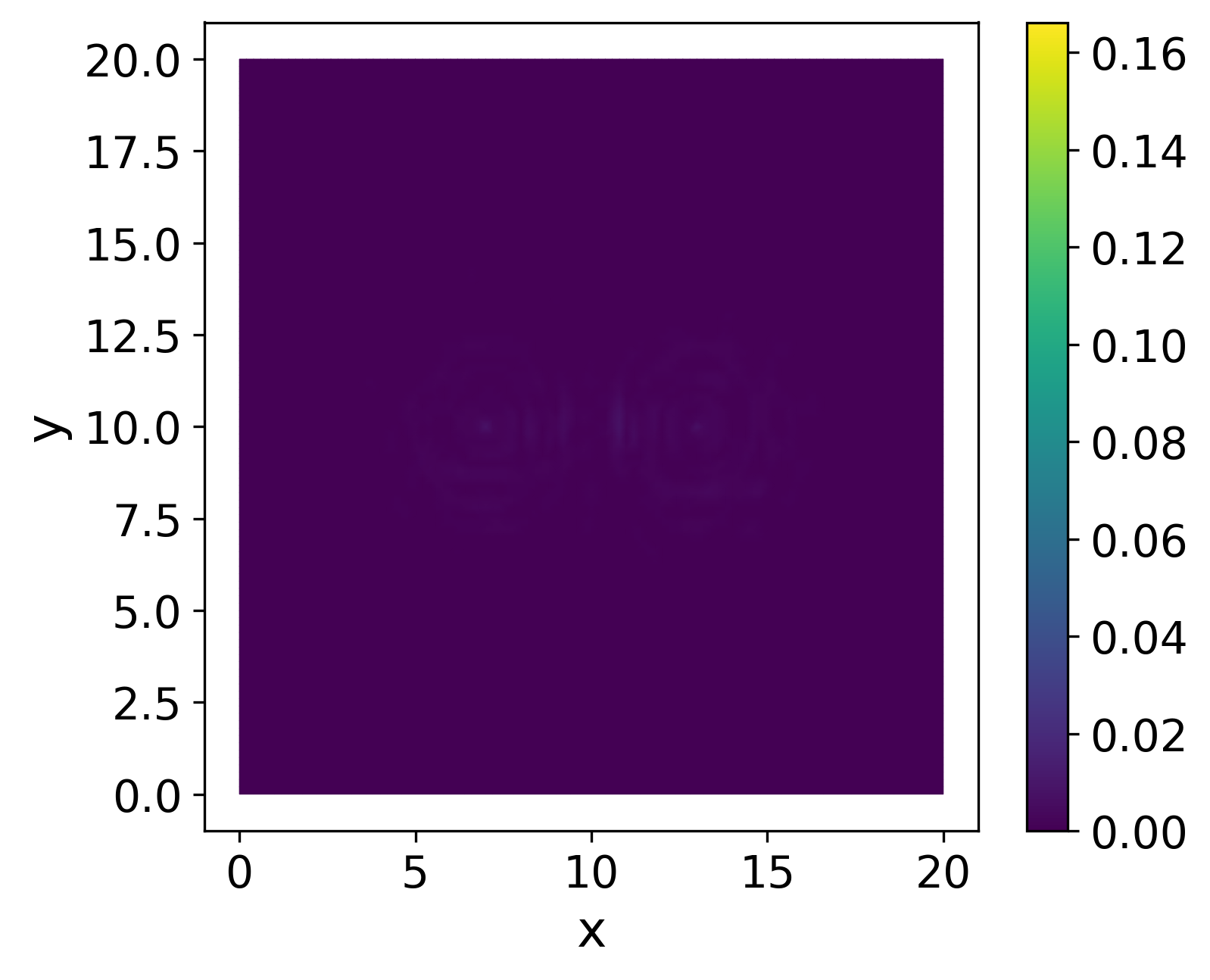}
    \caption{$t=0$.}
    \label{fig:test5-fields-t0}
  \end{subfigure}

  \vspace{0.4cm}

  \begin{subfigure}{\textwidth}
    \centering
    \includegraphics[width=0.32\textwidth]{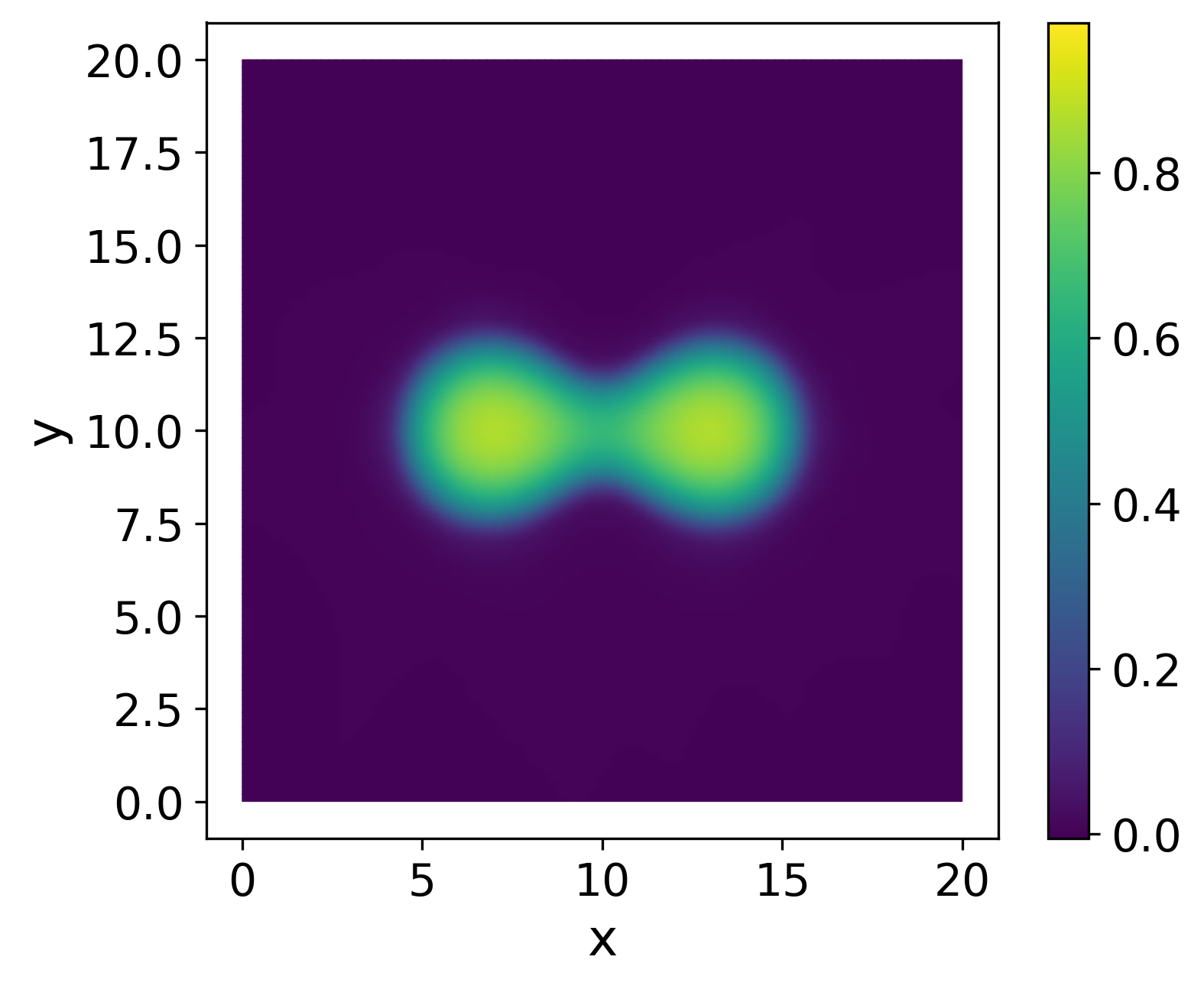}
    \hfill
    \includegraphics[width=0.32\textwidth]{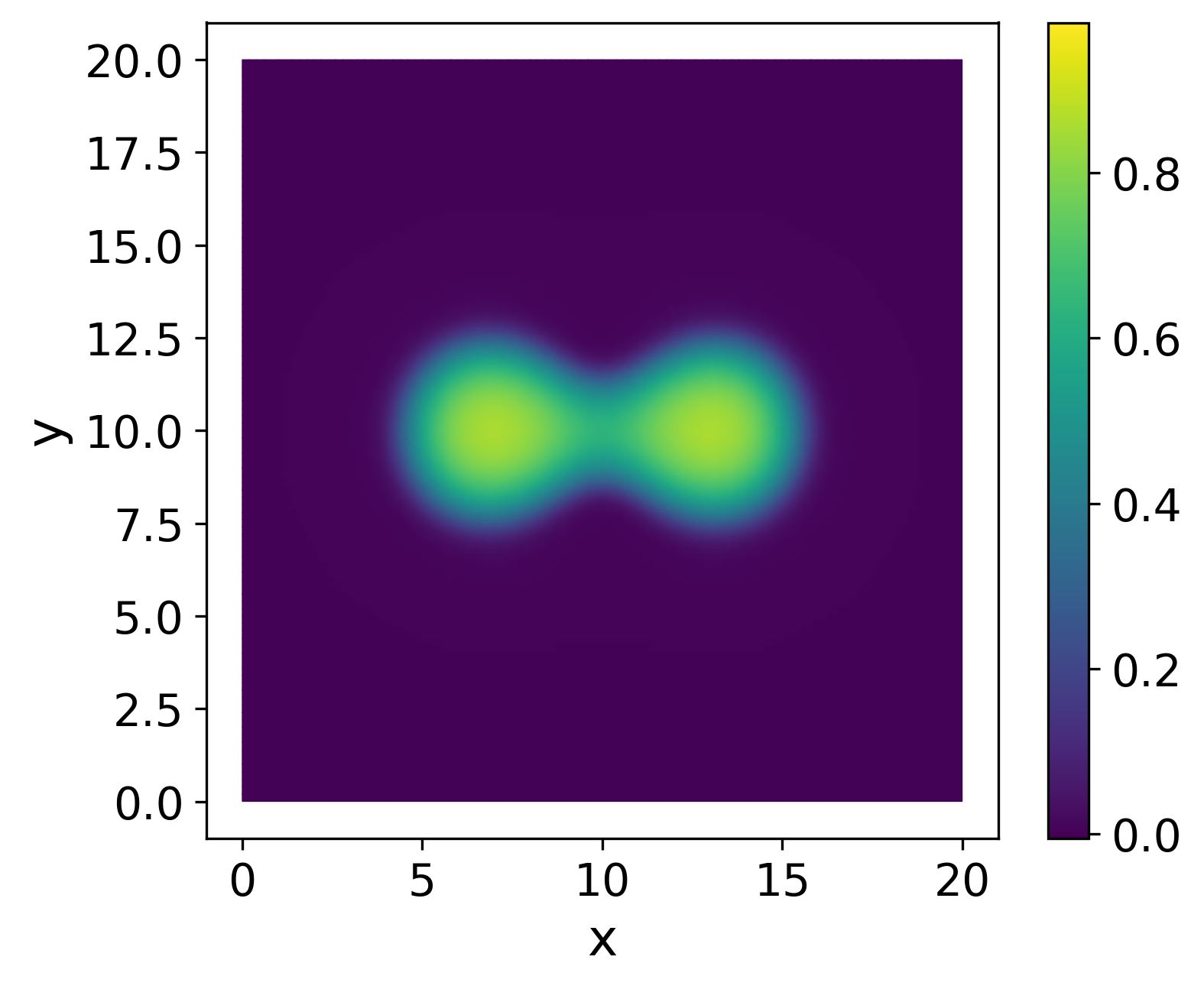}
    \hfill
    \includegraphics[width=0.32\textwidth]{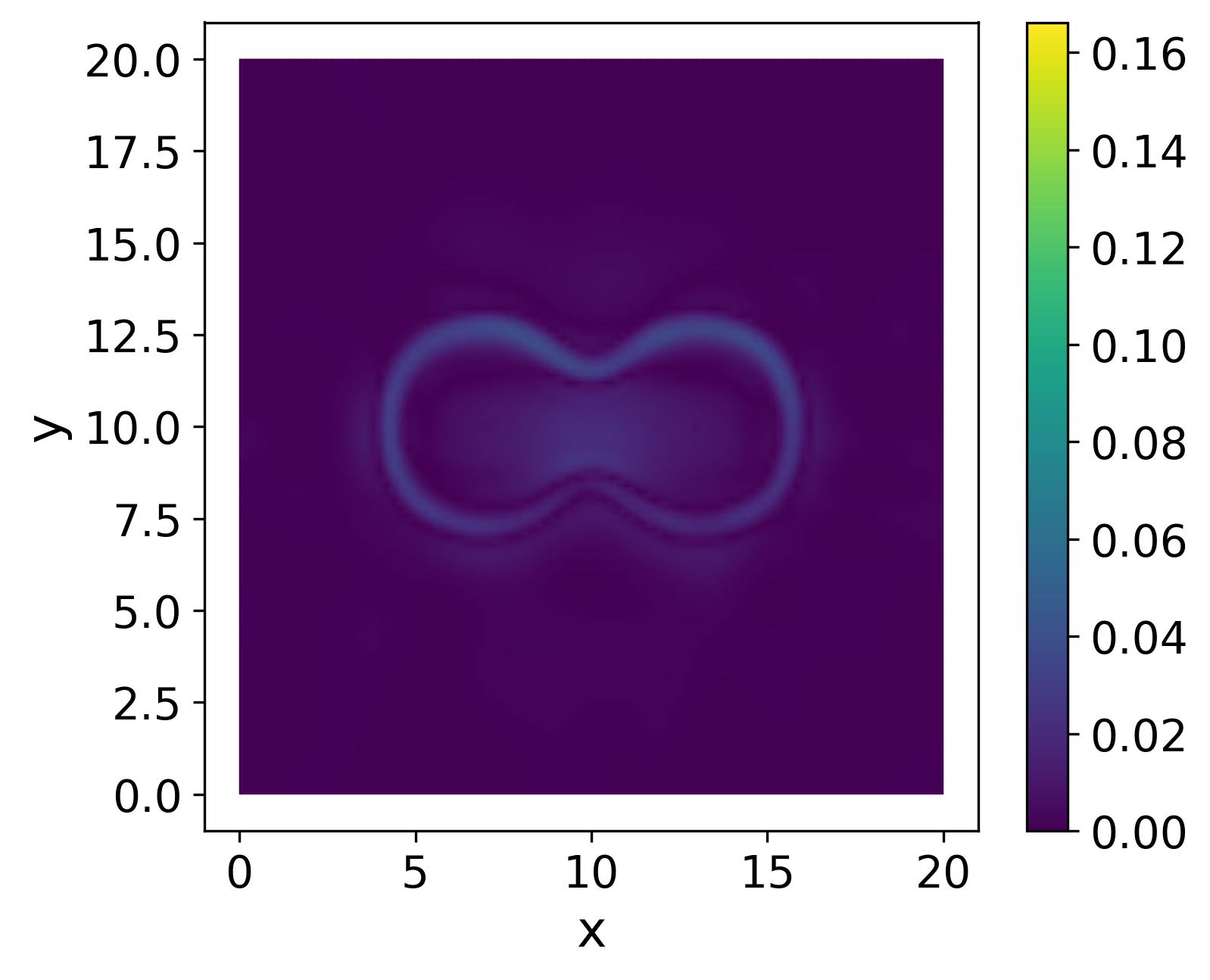}
    \caption{$t=1$.}
    \label{fig:test5-fields-t1}
  \end{subfigure}

  \vspace{0.4cm}

  \begin{subfigure}{\textwidth}
    \centering
    \includegraphics[width=0.32\textwidth]{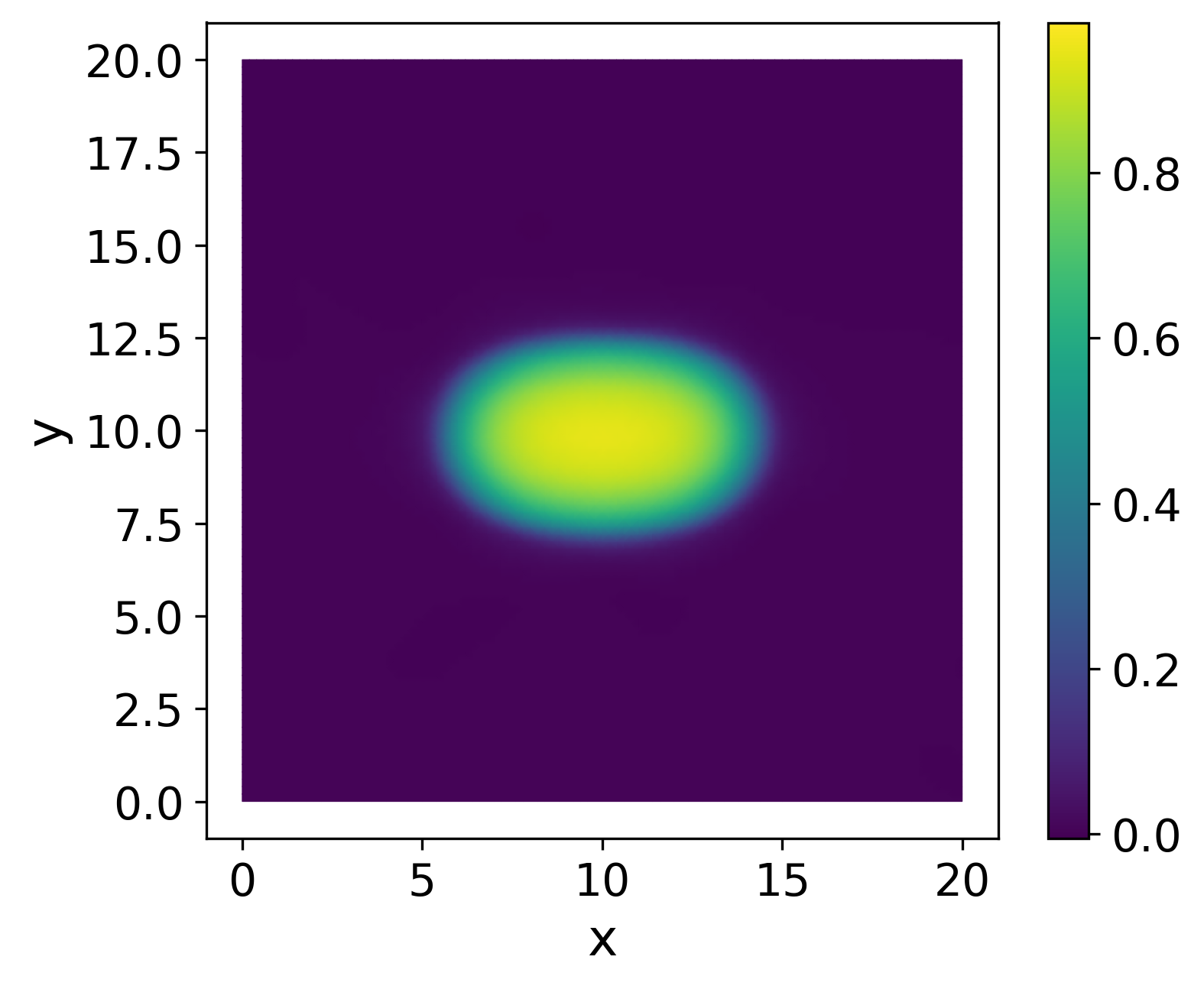}
    \hfill
    \includegraphics[width=0.32\textwidth]{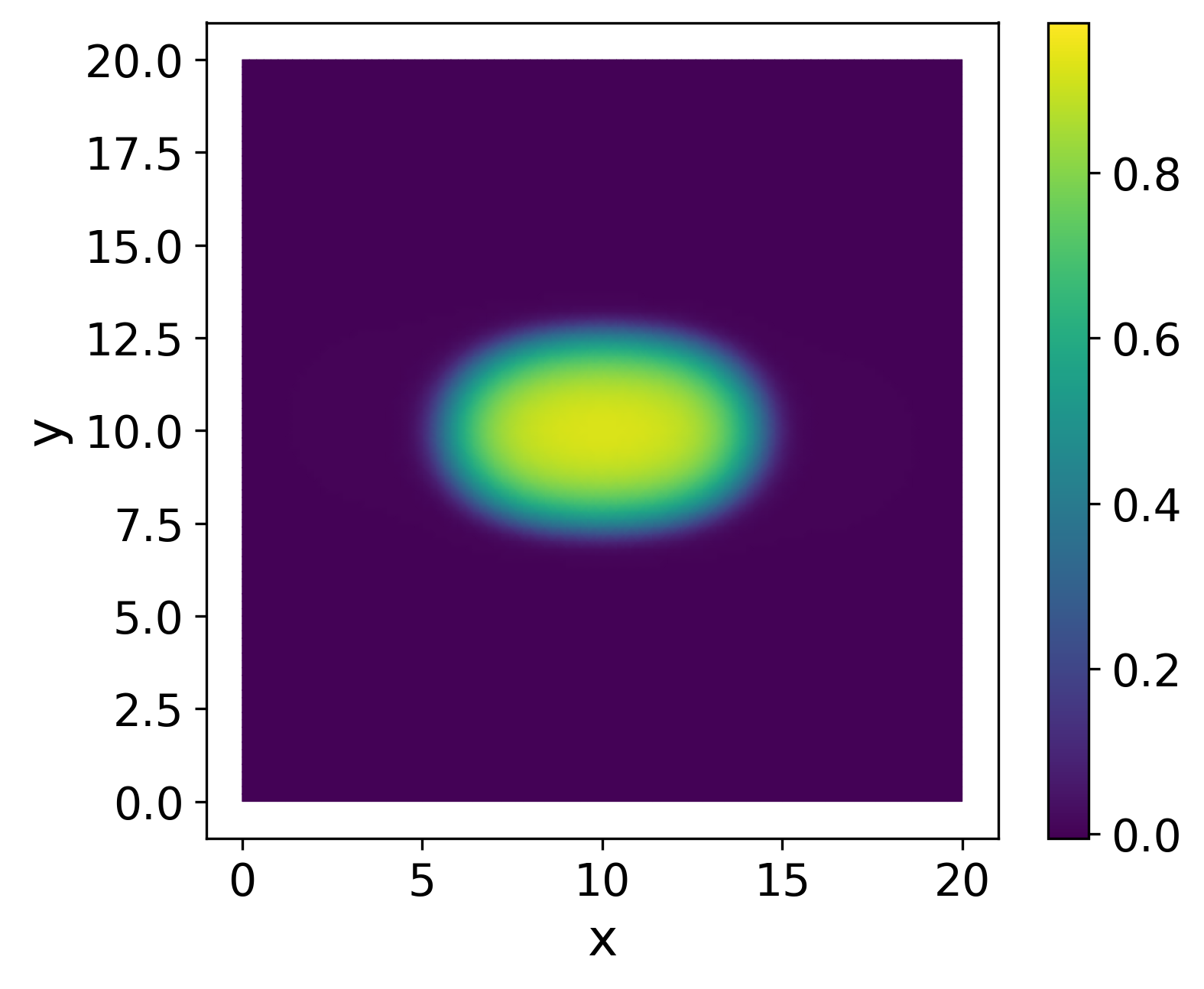}
    \hfill
    \includegraphics[width=0.32\textwidth]{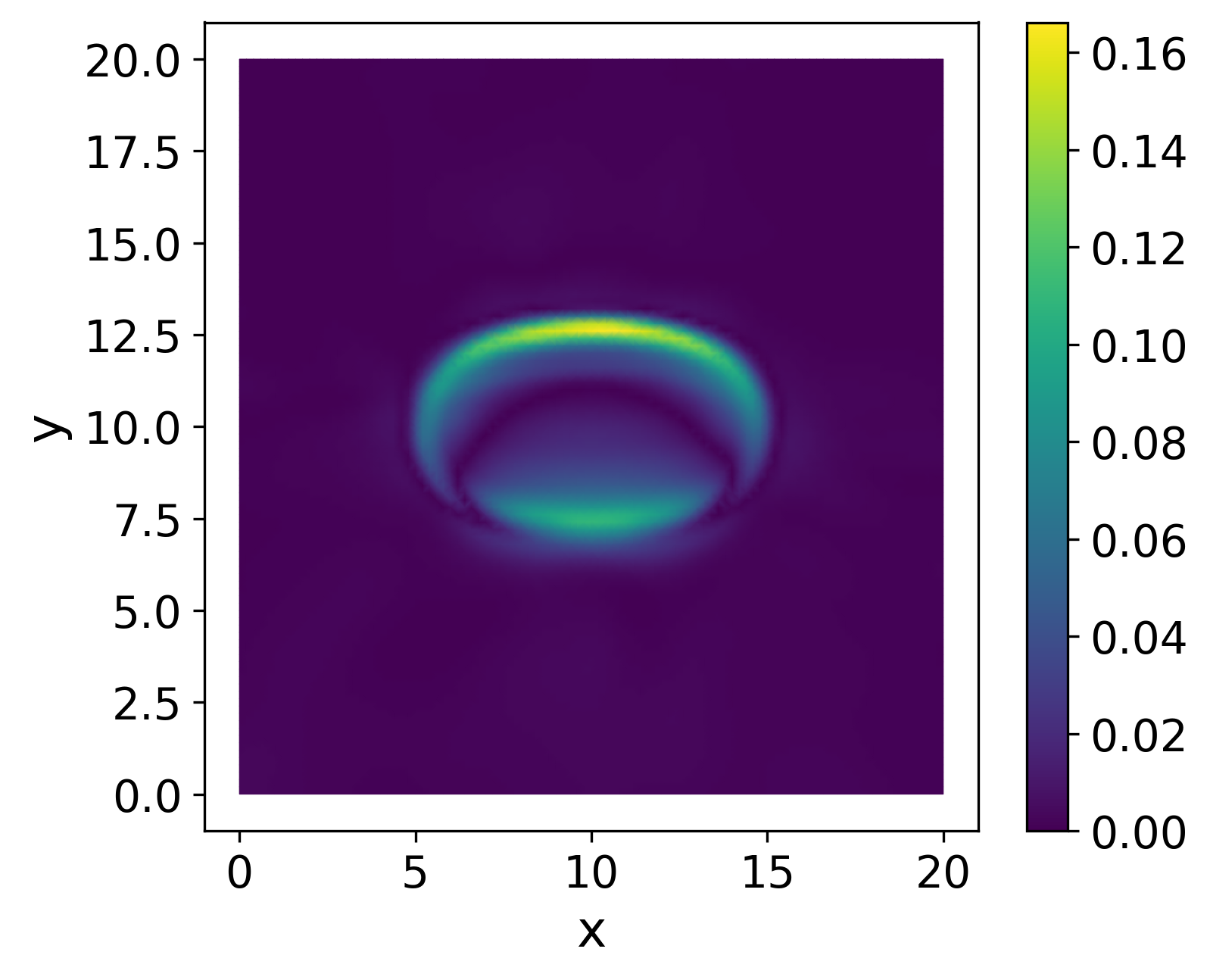}
    \caption{$t=7$.}
    \label{fig:test5-fields-t7}
  \end{subfigure}

  \caption{Two-dimensional density fields for the droplet-coalescence problem at three selected times. From left to right, the columns show the PINN approximation $\rho_\theta$, the CG-FEM reference solution $\rho_{\mathrm{ref}}$, and the absolute error $|\rho_\theta-\rho_{\mathrm{ref}}|$.}
  \label{fig:test5-fields-all}
\end{figure}

Figure~\ref{fig:test5-fields-all} compares the two-dimensional density fields obtained with the PINN and CG-FEM solutions, together with the corresponding absolute-error fields, at the three times considered in panel~(b) of Figure~\ref{fig:test5-2d-summary}.  At $t=0$, the PINN reproduces the initial two-droplet configuration accurately, with only very small localized discrepancies near the diffuse interfaces. At $t=1$, both methods capture the onset of coalescence and the formation of the bridge between the droplets, while the error becomes concentrated along the interfacial region. At $t=7$, both solutions have relaxed toward a single aggregated droplet, and the remaining discrepancy is again primarily localized near the interface.

\begin{figure}[h]
  \centering

  \begin{subfigure}{0.32\textwidth}
    \centering
    \includegraphics[width=\textwidth]{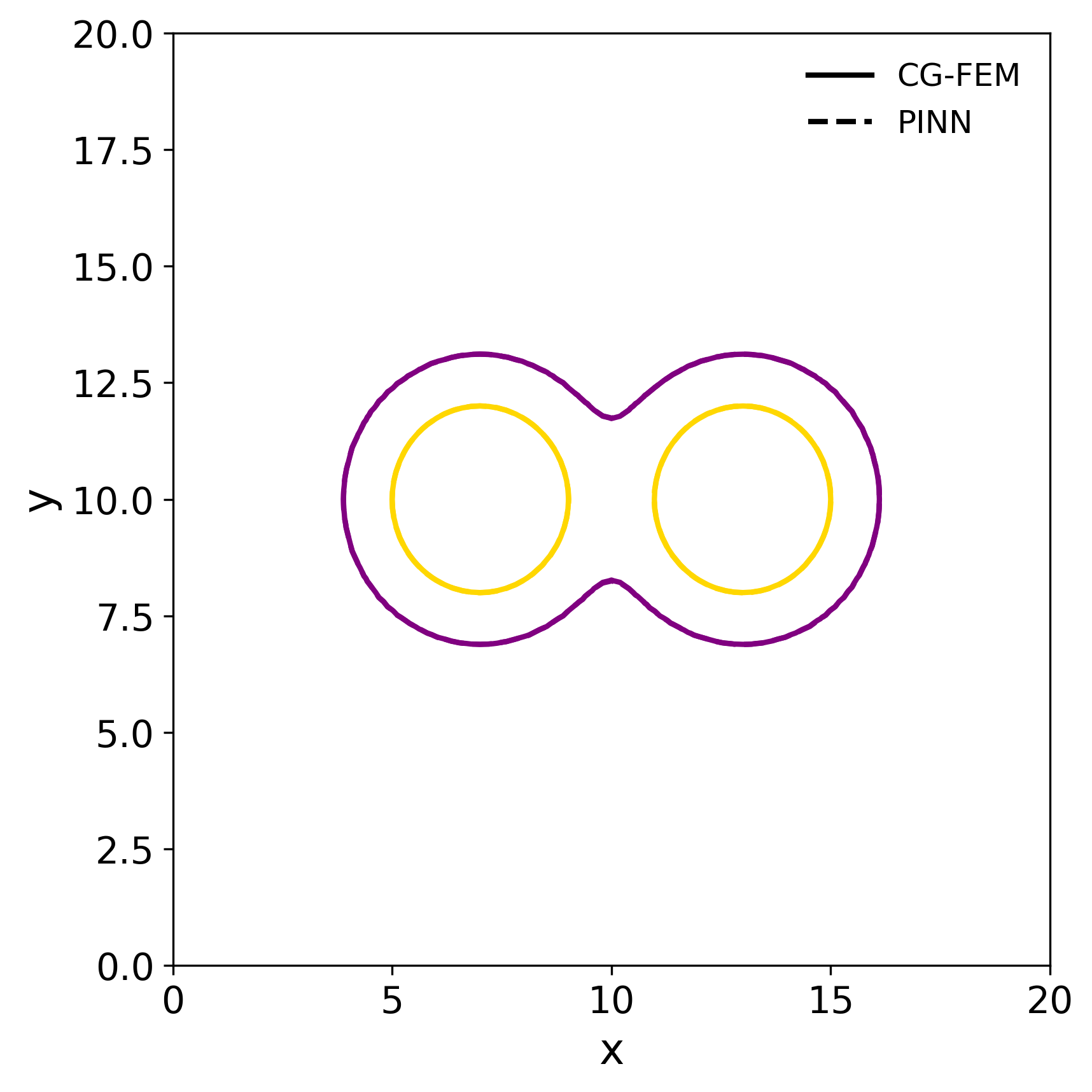}
    \caption{$t=0$.}
    \label{fig:test5-contour-t0}
  \end{subfigure}
  \hfill
  \begin{subfigure}{0.32\textwidth}
    \centering
    \includegraphics[width=\textwidth]{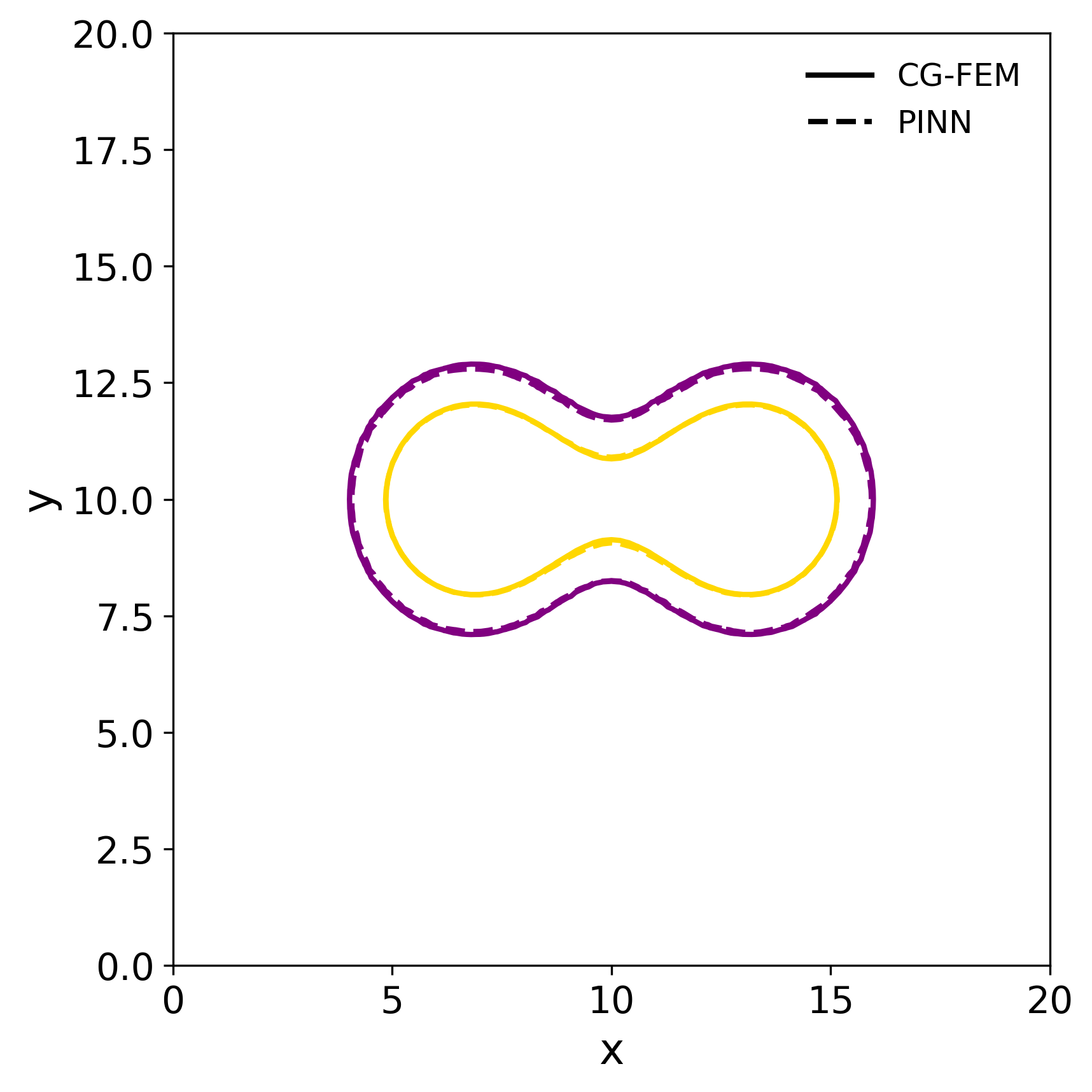}
    \caption{$t=1$.}
    \label{fig:test5-contour-t1}
  \end{subfigure}
    \hfill
  \begin{subfigure}{0.32\textwidth}
    \centering
    \includegraphics[width=\textwidth]{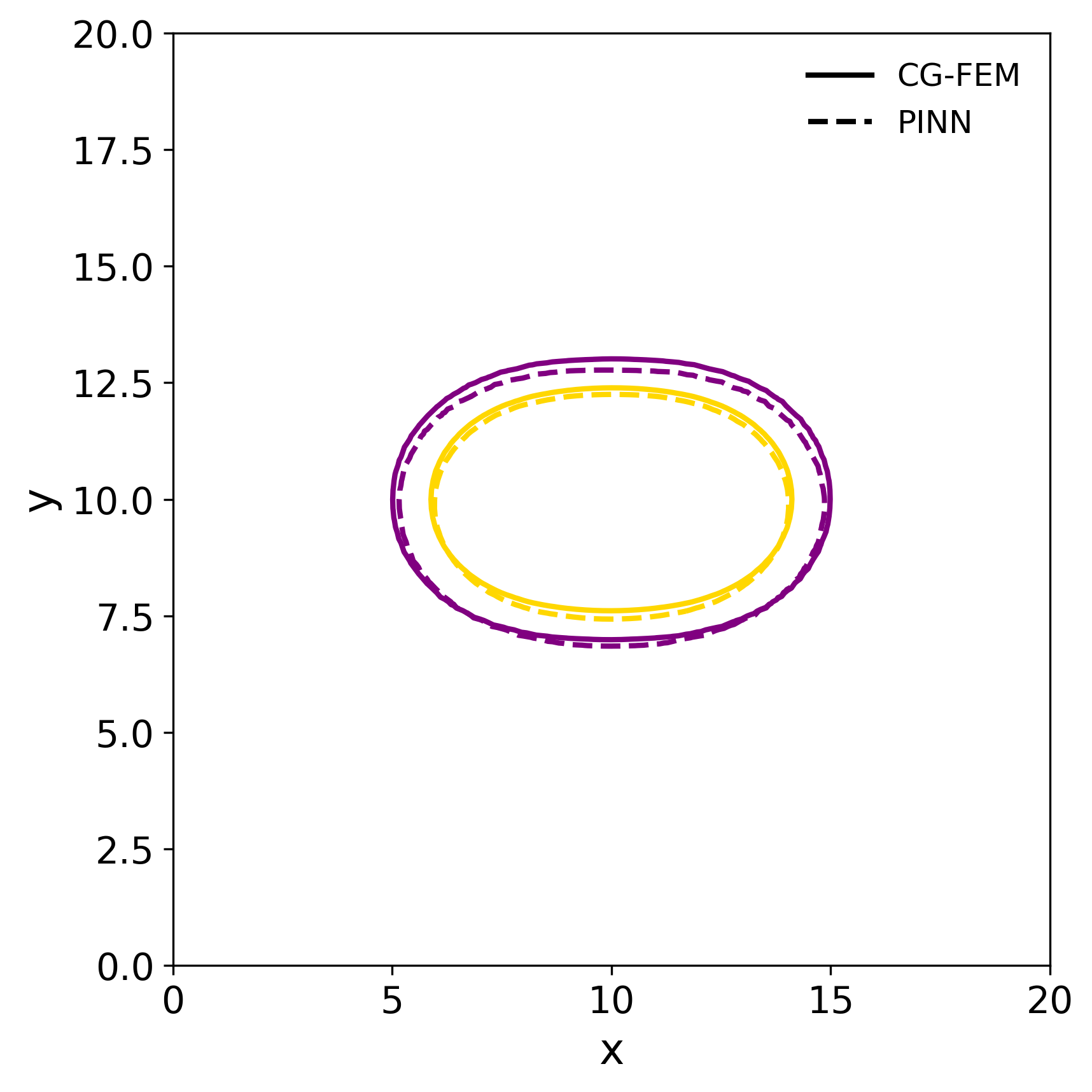}
    \caption{$t=7$.}
    \label{fig:test5-contour-t7}
  \end{subfigure}

\caption{Comparison of the density contours of the CG-FEM reference solution $\rho_{\mathrm{ref}}$ and the PINN approximation $\rho_\theta$ at $t=0$, $1$, and $7$. Solid lines represent $\rho_{\mathrm{ref}}$, while dashed lines represent $\rho_\theta$. The contours correspond to the density levels $\rho=0.1$ and $\rho=0.5$.}
  \label{fig:test5-contours-all}
\end{figure}

The contour plots in Figure~\ref{fig:test5-contours-all} confirm the same observations, namely that the agreement is already strong at the initial time, that the largest discrepancies occur during the transient coalescence stage near the bridge region, and that the two solutions come back into closer agreement once the merged droplet approaches equilibrium.

\section{Conclusions}\label{sec:conclusions}

In this work, we developed a physics-informed neural network framework for nonlocal partial differential equations arising in dynamic density functional theory. The study makes three principal contributions. First, it applies the PINN methodology to DDFT-type gradient-flow equations involving nonlinear diffusion, nonlocal interactions, mass conservation, and free-energy dissipation. Second, it introduces a modified Lorentzian activation function designed to improve training convergence and approximation accuracy. Third, it evaluates the nonlocal convolution term using a precomputed discrete operator, allowing the interaction contribution to be incorporated efficiently into the PDE residual during training.

The proposed framework was tested on one- and two-dimensional problems involving nonlocal interactions, nonlinear diffusion in an external potential, aggregation with attractive interactions, and droplet coalescence. The neural-network approximations were compared with CG-FEM and DG-FEM reference solutions and showed good agreement in terms of solution accuracy, error norms, and qualitative evolution. In comparison with the standard \texttt{tanh} activation function, the modified Lorentzian activation reached the target loss in fewer optimization steps and reduced the corresponding training time. The numerical results also reproduced the expected free-energy decay and preserved mass approximately, with only small deviations over the simulated time intervals.

The present study also has some limitations. First, the proposed PINN approach remains computationally more expensive than conventional FEM solvers for the forward problems considered here, particularly because training requires many optimization iterations and repeated evaluations of the PDE residual. Second, obtaining accurate solutions requires the calibration of several parameters, including the network architecture, loss weights, learning rate, activation parameter, batch size, and number of training iterations. The performance can therefore be sensitive to these choices, and a single parameter configuration may not work equally well across all examples. Finally, the current implementation of the precomputed nonlocal operator relies on a structured Cartesian grid. Consequently, the spatial nodes cannot be chosen arbitrarily, which limits the direct use of irregular collocation-point distributions.

These results indicate that PINNs can be applied to nonlocal DDFT equations when the network architecture and the evaluation of the interaction term account for the structure of the problem. The modified Lorentzian activation function and the precomputed convolution operator may also be useful for other nonlinear and nonlocal PDEs with similar approximation and optimization difficulties. Further work will consider longer time intervals, more strongly nonlocal regimes, and a more direct enforcement of mass conservation and free-energy dissipation. Additional extensions include weak-form PINN formulations and applications to more complex DDFT models and related gradient-flow equations.

\ack{The authors thank Peter Yatsyshin for valuable discussions and for the suggestion to represent the convolution term by a fixed matrix, as described in Appendix \ref{AppendixA}.}

\funding{

Dimitrios Gourzoulidis and Serafim Kalliadasis acknowledge financial support from the ERC--EPSRC Frontier Research Guarantee through Grant No.\ 588 EP/X038645, the ERC through Advanced Grant No.\ 247031, and the EPSRC through Grant Nos.\ EP/L025159 and EP/L020564.
\noindent
Soumaya Elkantassi acknowledges the following financial support: This material is based upon work supported by the Air Force Office of Scientific Research under award number FA8655-24-1-7009.

}

\roles{All authors contributed equally to this work.}



\appendix
\section{Appendix: discrete evaluation of the convolution term}\label{AppendixA}
This appendix explains how the nonlocal convolution term in \eqref{eq:variation} is evaluated in the finite-element discretization and how it is reduced to a discrete matrix operation in the implementation.

Let $V_h$ denote the finite-element
space, with basis functions $\{\phi_i\}_{i=1}^{N}$. We approximate the density
$\rho$ by $\rho_h\in V_h$, written as
\begin{equation*}
\rho_h(x)
=
\sum_{i=1}^{N}\rho_i\phi_i(x),
\end{equation*}
where $\rho_i$ are the coefficients of the finite-element approximation.

For a given interaction kernel $W$, the convolution of $W$ with $\rho_h$ is
defined by
\begin{equation*}
(W*\rho_h)(y)
=
\int_{\Omega}W(y-x)\rho_h(x)\,dx.
\end{equation*}
Evaluating this expression at a node $x_k$ gives
\begin{equation*}
(W * \rho_h)(x_k)
=
\sum_{i=1}^N \rho_i \int_\Omega W(x_k-x)\phi_i(x)\,dx.
\end{equation*}
To obtain a fully discrete representation, we approximate the shifted kernel $W(x_k-x)$ in the same finite-element basis by
\begin{equation*}
W(x_k-x)\approx \sum_{m=1}^N B_{mk}\phi_m(x),
\end{equation*}
where $B_{mk}$ are the corresponding interpolation or projection coefficients. Substituting this approximation into the convolution expression yields
\begin{equation*}
(W * \rho_h)(x_k)
\approx
\sum_{i=1}^N \sum_{m=1}^N
\rho_i B_{mk}\int_\Omega \phi_i(x)\phi_m(x)\,dx.
\end{equation*}
Introducing the mass matrix
\begin{equation*}
M_{im}=\int_\Omega \phi_i(x)\phi_m(x)\,dx,
\end{equation*}
we obtain
\begin{equation*}
(W * \rho_h)(x_k)
\approx
\sum_{i=1}^N \sum_{m=1}^N \rho_i M_{im} B_{mk}.
\end{equation*}
Since the mass matrix is symmetric, \(M_{im}=M_{mi}\), this becomes
\begin{equation*}
(W * \rho_h)(x_k)
\approx
\sum_{i=1}^N \sum_{m=1}^N B_{mk} M_{mi}\rho_i
=
(B^\T M \boldsymbol{\rho})_k,
\end{equation*}
where \(\boldsymbol{\rho}=(\rho_1,\dots,\rho_N)^\T\). Hence, with
\begin{equation*}
\Phi = B^\T M,
\end{equation*}
the discrete convolution is evaluated through the matrix--vector product \(\Phi\boldsymbol{\rho}\).

\section{Appendix: reference finite element solvers}\label{AppendixB}

This appendix summarizes the finite-element reference solvers used in the numerical experiments of Section~\ref{sec:num} to generate the solutions against which the neural-network approximations are compared. For completeness, we first recall the weak formulation of \eqref{eq:equivalent}, which is the standard variational form used as the starting point for finite-element discretization; see, for example, \cite{quarteroni2008numerical,ernguermond2004theory}.

Starting from \eqref{eq:equivalent}, we multiply by a test function $v$ and integrate over $\Omega$ to obtain the weak formulation
\begin{equation*}
\int_\Omega \partial_t \rho \, v \, dx
+ \int_\Omega \rho \,\nabla \mu \cdot \nabla v \, dx = 0
\qquad \forall v \in V,
\end{equation*}
where $V$ is a suitable test space. The boundary contribution vanishes under the no-flux condition $\rho \,\nabla \mu \cdot \mathbf{n}=0$ on $\partial\Omega$.

\subsection{Continuous Galerkin formulation}

We first consider a continuous Galerkin discretization, which approximates the solution by functions that remain continuous across element interfaces. This is the reference method used for the two-dimensional computations in Section~\ref{sec:num}.

Let $V_h \subset H^1(\Omega)$ be a continuous finite-element space, for example the space of continuous piecewise linear functions on a conforming mesh of $\Omega$. Using a semi-implicit Crank--Nicolson discretization in time, the fully discrete problem is: given $\rho_h^n \in V_h$ at time $t^n$, find $\rho_h^{n+1}\in V_h$ such that
\begin{equation*}
\int_\Omega \dfrac{\rho_h^{n+1}-\rho_h^n}{\Delta t}\, v_h\, dx
+ \int_\Omega \rho_h^{n+1/2}\,\nabla \mu_h^{n+1/2}\cdot \nabla v_h\, dx = 0
\qquad \forall v_h \in V_h,
\end{equation*}
where
\begin{equation*}
\rho_h^{n+1/2}=\frac12\left(\rho_h^{n+1}+\rho_h^n\right),
\qquad
\mu_h^{n+1/2}=H'\!\left(\rho_h^{n+1/2}\right)+V+W * \rho_h^{n}.
\end{equation*}
Here the convolution term is evaluated explicitly at the previous time step,
\begin{equation*}
(W * \rho_h^{n})(\mathbf{x})
=
\int_\Omega W(\mathbf{x}-\mathbf{y})\,\rho_h^{n}(\mathbf{y})\,d\mathbf{y}.
\end{equation*}
Thus, the local nonlinear terms are treated at the midpoint level, whereas the nonlocal interaction term is treated explicitly. The nonlinear system arising at each time step is solved by Newton iteration, using the Fr\'echet derivative with respect to $\rho_h^{n+1}$; see, for example, \cite{CaboussatGourzoulidis2021}.

\subsection{Discontinuous Galerkin formulation}

For the one-dimensional test cases of Section~\ref{sec:num}, we also consider a discontinuous Galerkin discretization. In this approach, the approximate solution is allowed to be discontinuous across element interfaces, and communication between neighboring elements is enforced through numerical fluxes and jump terms.

Let $X_h$ be a discontinuous finite-element space of piecewise polynomial functions on the mesh. Using the same semi-implicit Crank--Nicolson time discretization, the DG-FEM scheme is: given $\rho_h^n \in X_h$, find $\rho_h^{n+1}\in X_h$ such that
\begin{align*}
\int_\Omega \dfrac{\rho_h^{n+1}-\rho_h^n}{\Delta t}\,v_h\,dx
&+ \int_\Omega \rho_h^{n+1/2}\,\nabla \mu_h^{n+1/2}\cdot \nabla v_h\,dx \\
&- \sum_{e\in\mathcal{E}_{\mathrm{int}}}\int_e (\rho_h^{n+1/2}\nabla \mu_h^{n+1/2})^{*}\,[v_h]\,ds \\
&+ 2\sum_{e\in\mathcal{E}_{\mathrm{int}}}\int_e [\rho_h^{n+1/2}]\,[v_h]\,ds
=0,
\qquad \forall v_h\in X_h.
\end{align*}
Here $\mathcal{E}_{\mathrm{int}}$ denotes the set of interior interfaces of the mesh, that is, the boundaries shared by neighboring elements, $[q]=q^+-q^-$ is the jump of a scalar quantity across an interior interface, and
\begin{equation*}
\rho_h^{n+1/2}=\frac12\left(\rho_h^{n+1}+\rho_h^n\right),
\qquad
\mu_h^{n+1/2}=H'\!\left(\rho_h^{n+1/2}\right)+V+W * \rho_h^{n}.
\end{equation*}
In one space dimension, the upwind numerical flux is defined by
\begin{equation*}
(\rho_h^{n+1/2}\beta_h)^*
=
\max(\hat{\beta}_h,0)\,\rho_h^{n+1/2,-}
+
\min(\hat{\beta}_h,0)\,\rho_h^{n+1/2,+},
\qquad
\beta_h=\partial_x\mu_h^{n+1/2},
\end{equation*}
where $\hat{\beta}_h$ denotes the average of the traces of $\beta_h$ across the interface. The last term is a jump stabilization term. As in the continuous Galerkin case, the nonlinear system is solved by Newton iteration. After each time step, a minmod-type slope limiter is applied to reduce spurious oscillations; see \cite{hesthavenwarburton2007dg}.

\newpage
\bibliographystyle{iopart-num}
\bibliography{references}

\end{document}

%% file: macros.tex
\def\frac#1#2{{\textstyle{#1\over#2}}}

\def\bi{\begin{itemize}}
\def\ei{\end{itemize}}
\def\bd{\begin{description}}
\def\ed{\end{description}}
\def\ben{\begin{enumerate}}
\def\een{\end{enumerate}}
\def\bv{\begin{verbatim}}
\def\ev\end{verbatim}

\def\T{{ \mathrm{\scriptscriptstyle T} }}

\DeclareMathOperator*{\E}{E}

\def\2to{{\ {\buildrel 2\over \longrightarrow}\ }}

\def\I1ton{{$I_1,\ldots,I_n$}}
\def\X1ton{{$X_1,\ldots,X_n$}}
\def\Y1ton{{$Y_1,\ldots,Y_n$}}
\def\Z1ton{{$Z_1,\ldots,Z_n$}}
\def\R1ton{{$R_1,\ldots,R_n$}}
\def\e1ton{{$e_1,\ldots,e_n$}}
\def\t1ton{{$t_1,\ldots,t_n$}}
\def\x1ton{{$x_1,\ldots,x_n$}}
\def\y1ton{{$y_1,\ldots,y_n$}}
\def\z1ton{{$z_1,\ldots,z_n$}}